\documentclass[12pt]{article}
\usepackage{a4wide,epsf,epsfig,amssymb,amsmath}
\usepackage{bbm,theorem}
\newcommand{\one}{\mathbf{1}}
\newcommand{\zero}{\mathbf{0}}
\newcommand{\aaa}{\mbox{{\boldmath$a$}}}
\newcommand{\bb}{\mbox{{\boldmath$b$}}}
\newcommand{\ee}{\mbox{{\boldmath$e$}}} 
\newcommand{\cc}{\mbox{{\boldmath$c$}}}
\newcommand{\hh}{\mbox{{\boldmath$h$}}} 
\newcommand{\uu}{\mbox{{\boldmath$u$}}} 
\newcommand{\vv}{\mbox{{\boldmath$v$}}} 
\newcommand{\ww}{\mbox{{\boldmath$w$}}} 
\newcommand{\xx}{\vphantom{l}{\mbox{{\boldmath$x$}}}} 
\newcommand{\subxx}{\mbox{{\boldmath$_x$}}}
\newcommand{\smallxx}{\mbox{{\boldmath\scriptsize $x$}}}
\newcommand{\supcc}{\mbox{{\boldmath$^c$}}}
\newcommand{\const}{\mbox{{\boldmath$\gamma$}}}
\newcommand{\yy}{\mbox{{\boldmath$y$}}} 
\newcommand{\N}{{\mathbbm N}}
\newcommand{\Z}{{\mathbbm Z}}
\newcommand{\Q}{{\mathbbm Q}}
\newcommand{\R}{{\mathbbm R}}

\newcommand{\Ahat}{{\widehat A}}

\newcommand{\symmdiff}{{\vartriangle}}
\newcommand{\conv}{\mbox{\rm conv}}
\newcommand{\diam}{\mbox{\rm diam}}
\newcommand{\aff}{\mbox{\rm aff}}
\newcommand{\Vol}{\mbox{\rm{Vol}}}
\newcommand{\CUT}{\mbox{\rm{CUT}}}
\newcommand{\MET}{\mbox{\rm{MET}}}
\newcommand{\COR}{\mbox{\rm{COR}}}
\newcommand{\query}[1]{\marginpar{?? #1}}
\theoremstyle{break}%
{\theorembodyfont{\itshape}%
\newtheorem{Theorem}{Theorem}%
\newtheorem{Proposition}[Theorem]{Proposition}%
\newtheorem{Corollary}[Theorem]{Corollary}
\newtheorem{Lemma}[Theorem]{Lemma}%
\newtheorem{Conjecture}[Theorem]{Conjecture}}%
{\theorembodyfont{\rmfamily}%
\newtheorem{Definition}[Theorem]{Definition}%
\newtheorem{Example}[Theorem]{Example}%
\newtheorem{Notation}[Theorem]{Notation}%
\newtheorem{ExampleExercise}[Theorem]{Example/Exercise}%
}%

\newcommand\proof{\par\medbreak\noindent{\bf Proof.}\enspace}
\newcommand\qed{\vbox{\hrule
  \hbox{\vrule\hbox to 5pt{\vbox to 8pt{\vfil}\hfil}\vrule}\hrule}}
\newcommand\proofend{\unskip \nobreak \hskip0pt plus 1fill \qquad \qed
\medskip\noindent}
\newcommand\sm{{\setminus}}
\newcommand\sse{\subseteq}
\begin{document}
\title{\bfseries Lectures on 0/1-Polytopes}
\author{{\LARGE G\"unter M. Ziegler\thanks{%
Supported by a DFG Gerhard-Hess-Forschungsf\"orderungspreis (Zi 475/1-2)
and by a German Israeli Foundation (G.I.F.) grant
I-0309-146.06/93.}}\\[1mm]
Dept.\ Mathematics, MA~7-1\\ 
Technische Universit\"at Berlin\\
10623 Berlin, Germany\\
{\tt ziegler@math.tu-berlin.de}}
\maketitle

\begin{abstract}
These lectures on the combinatorics and geometry of
$0/1$-polytopes are meant as an \emph{introduction\/} and
\emph{invitation}. Rather than heading for an extensive survey on
$0/1$-polytopes I present some interesting aspects of these
objects; all of them are related to some quite recent work and
progress.

$0/1$-polytopes have a very
simple definition and explicit descriptions;
we can enumerate and analyze small
examples explicitly in the computer (e.~g.\ using {\tt polymake}).
However, any intuition that is derived from the
analysis of examples in ``low dimensions''
will miss the true complexity of $0/1$-polytopes.
Thus, in the following we will study several aspects 
of the complexity of higher-dimensional $0/1$-polytopes:
the doubly-exponential number of combinatorial types, 
the number of facets which can be huge, and the coefficients of
defining inequalities which sometimes turn out to be extremely large.
Some of the effects and results will be backed by proofs
in the course of these lectures; we will also 
be able to verify some of them on
explicit examples, which are accessible as a {\tt polymake}
database.
\end{abstract}

\newpage
\begin{small}\tableofcontents\end{small}

\section*{Introduction}
\addcontentsline{toc}{section}{Introduction}

These lectures are trying to get you interested in
$0/1$-polytopes. But I must warn you: they are mostly
``bad news lectures'' --- with 
two types of bad news:
\begin{enumerate}
\item
General $0/1$-polytopes are complicated objects, and some of them have various
kinds of extremely bad properties such as 
``huge coefficients'' and ``many facets,''
which are bad news also with respect to applications.
\item 
Even worse, there are bad gaps in our understanding of $0/1$-polytopes.
Very basic problems and questions are open,
some of them embarassingly 
easy to state, but hard to answer.
So, $0/1$-polytopes are interesting and remain {\em challenging}.
\end{enumerate}

A good grasp on the structure of  
$0/1$-polytopes is important for the ``polyhedral combinatorics''
approach of combinatorial optimization.
This has motivated an extremely thorough study of 
some special classes of $0/1$-polytopes such as the traveling
salesman polytopes (see Gr\"otschel \& Padberg \cite{GroetschelPadberg}
and Applegate, Bixby, Cook \& Chv\'atal \cite{ABCC2})
and the cut polytopes (see Deza \& Laurent \cite{DL}, and 
Section~\ref{sec:cut}).
In such studies the question about properties of
general $0/1$-polytopes, and for complexity estimates
about them, arises quite frequently and naturally.
Thus Gr\"otschel \& Padberg \cite{GroetschelPadberg}
looked for upper bounds on the number of facets,
and we can now considerably improve the estimates they 
had then (Section~\ref{sec:facets}).
One also asks for the sizes of the integers that appear as
facet coefficients --- and the fact that these coefficients
may be huge (Section~\ref{sec:coeff}) is bad news 
since it means that there is a great danger of numerical instability
or arithmetic overflow.

Surprisingly, however, properties of {\em general $0/1$-polytopes\/}
have not yet been a focus of research.
I think they should be, and these lecture notes (expanded from my
DMV-Seminar lectures in Oberwolfach, November 1997)
are meant to provide support for this.

Of course, the distinction between ``special'' and ``general''
$0/1$-polytopes is somewhat artificial.
For example, Billera \& Sarangarajan \cite{BS} have 
proved the surprising fact
that {\em every\/} $0/1$-polytope appears as a face of a
TSP-polytope.
Nevertheless, a study of the broad class of general
$0/1$-polytopes provides new points of view. Here it
appears natural to look at {\em extremal\/} polytopes 
(e.~g.\ polytopes with ``many facets''), and at
{\em random polytopes\/} and their properties.

Where is the difficulty in this study? 
The definition of $0/1$-polytopes is very simple,
examples are easy to come by, and they can be analyzed completely.
But this simplicity is misleading:
there are various effects that appear only in rather high 
dimensions ($d\gg3$, whatever that means). Part of this we
will trace to one basic linear algebra concept:
determinants of $0/1$-matrices, which show their typical
behaviour --- large values, and a low probability to vanish ---
only when the dimension gets quite large.
Thus one rule of thumb will be justified again and again:
\begin{center}
{\sf Low-dimensional intuition does not work!}
\end{center}
Despite this (and to demonstrate this), our discussion in
various lectures will take the low-dimensional situation
as a starting point, and as a point of reference.
(For example, the first lecture will start with a
list of $3$-dimensional $0/1$-polytopes,
which will turn out to be deceptively simple.)

However, examples are nevertheless important.
The {\tt polymake} project \cite{GawrilowJoswig,GawrilowJoswig2}
provides a framework and many fundamental tools
for their detailed analysis.
Thus, these lecture notes come with a library of
interesting examples, provided as a separate section of the
{\tt polymake} database at
\begin{itemize}
\item[]
{\tt http://www.math.tu-berlin.de/diskregeom/polymake/}  
\end{itemize}
We will refer to examples in this database throughout.
The names of the polytope data files are of the form
{\tt NN:d-n.poly},
where {\tt NN} is an identifyer of the polytope
(e.~g.\ initials of whoever supplied the example),
{\tt d} is the dimension of the polytope, {\tt n} is its 
number of vertices.
I invite you to play with these examples.
(Also, I am happy to accept further contributions to extend this
bestiary of interesting $0/1$-polytopes!)
\section{Classification of Combinatorial Types}

\subsection{Low-dimensional $0/1$-polytopes}

$0/1$-polytopes may be defined as the convex hulls
of finite sets of $0/1$-vectors, that is, as the convex
hulls of subsets of the vertices of the regular cube
$C_d=\{0,1\}^d$. Until further notice let's assume that
we only consider full-dimensional $0/1$-polytopes, so we
have $P=P(V)=\conv(V)$ for some $V\sse\{0,1\}^d$,
where we assume that $P$ has dimension~$d$.
We call two polytopes {\em $0/1$-equivalent\/} if one can be
transformed into the other by a symmetry of the $0/1$-cube.

Now $0/1$-polytopes of dimensions $d\le2$ are not interesting:
we get a point, the interval $[0,1]$, a triangle, and a square.
\[
\def\epsfsize#1#2{.5#1}
\epsffile{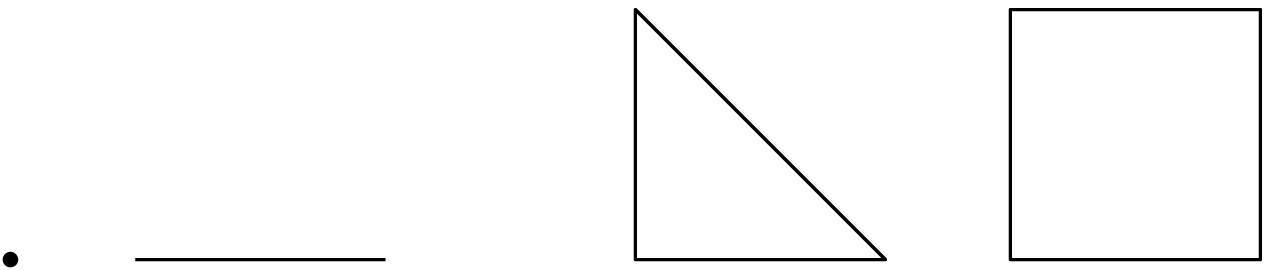}
\]
The figure below represents the classification of 
$3$-dimensional $0/1$-polytopes
$P\subseteq\R^3$ according to $0/1$-equivalence. An arrow
$P\epsffile{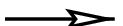} P'$
between two of them denotes that $P'$ is $0/1$-equivalent to a
\emph{subpolytope} of $P$, that is, $P'\sim P(V')$ and $P=P(V)$ for
some subset $V'\subseteq V$.

\newpage
~~\vskip-1.6cm
\[
\begin{picture}(0,0)%
\epsfig{file=EPS/3Class2.pstex}%
\end{picture}%
\setlength{\unitlength}{2763sp}%
\begingroup\makeatletter\ifx\SetFigFont\undefined%
\gdef\SetFigFont#1#2#3#4#5{%
  \reset@font\fontsize{#1}{#2pt}%
  \fontfamily{#3}\fontseries{#4}\fontshape{#5}%
  \selectfont}%
\fi\endgroup%
\begin{picture}(9783,9547)(-981,-9136)
\put(301,-3961){\makebox(0,0)[lb]{\smash{\SetFigFont{12}{14.4}{\rmdefault}{\mddefault}{\updefault}prism}}}
\put(3601,-9136){\makebox(0,0)[lb]{\smash{\SetFigFont{12}{14.4}{\rmdefault}{\mddefault}{\updefault}COR(2)}}}
\put(7801,-9136){\makebox(0,0)[lb]{\smash{\SetFigFont{12}{14.4}{\rmdefault}{\mddefault}{\updefault}CUT(3)}}}
\put(6001,-9136){\makebox(0,0)[lb]{\smash{\SetFigFont{12}{14.4}{\rmdefault}{\mddefault}{\updefault}$\Delta'_3$}}}
\put(2551,-5911){\makebox(0,0)[lb]{\smash{\SetFigFont{12}{14.4}{\rmdefault}{\mddefault}{\updefault}another}}}
\put(2551,-6136){\makebox(0,0)[lb]{\smash{\SetFigFont{12}{14.4}{\rmdefault}{\mddefault}{\updefault}square}}}
\put(2551,-6361){\makebox(0,0)[lb]{\smash{\SetFigFont{12}{14.4}{\rmdefault}{\mddefault}{\updefault}pyramid}}}
\put(2476,-3961){\makebox(0,0)[lb]{\smash{\SetFigFont{12}{14.4}{\rmdefault}{\mddefault}{\updefault}nameless}}}
\put(2926,-361){\makebox(0,0)[lb]{\smash{\SetFigFont{12}{14.4}{\rmdefault}{\mddefault}{\updefault}cube}}}
\put(7051,-6361){\makebox(0,0)[lb]{\smash{\SetFigFont{12}{14.4}{\rmdefault}{\mddefault}{\updefault}bipyramid}}}
\put(7276,-3961){\makebox(0,0)[lb]{\smash{\SetFigFont{12}{14.4}{\rmdefault}{\mddefault}{\updefault}octahedron}}}
\put(151,-6361){\makebox(0,0)[lb]{\smash{\SetFigFont{12}{14.4}{\rmdefault}{\mddefault}{\updefault}pyramid}}}
\put(151,-6136){\makebox(0,0)[lb]{\smash{\SetFigFont{12}{14.4}{\rmdefault}{\mddefault}{\updefault}square}}}
\put(526,-8536){\makebox(0,0)[rb]{\smash{\SetFigFont{12}{14.4}{\rmdefault}{\mddefault}{\updefault}$4$ tetrahedra:}}}
\end{picture}

\]

The full-dimensional
$0/1$-polytopes of dimension $d=4$ were first enumerated by Alexx
Below: There are $349$ different $0/1$-equivalence classes.

In dimension $5$ there are exactly
$1226525$ different $0/1$-equivalence classes
of $5$-dimensional $0/1$-polytopes. This classification
was done by Oswin Aichholzer \cite{Ai2}:
a considerable achievement, which was possible only by 
systematic use of all the symmetry that is inherent in the problem.

In October 1998, Aichholzer completed also an enumeration
and classification of the $6$-\allowbreak di\-men\-sional $0/1$-polytopes
up to $12$ vertices. 
%
The complete classification of all $6$-dimensional $0/1$-polytopes is not
within reach: in fact, even the output, a non-redundant list
of all combinatorial types
would be so huge that it is impossible to store or search
efficiently: and thus it would probably\footnote{%
``Where a calculator like the ENIAC today is equipped with $18,000$ vacuum
tubes and weighs 30 tons, computers in the future may have only $1,000$
vacuum tubes and perhaps weigh only $1\frac12$ tons.''
--- {\sl Popular Mechanics}, March 1949, p.~258.}
be useless.

\subsection{Combinatorial types}

Many fundamental concepts of general polytope theory 
can be specialized to the situation of $0/1$-polytopes.
The following reviews the basic definitions and concepts.
See for example \cite[Lect.~0-3]{Zpoly} for 
more detailed explanations.

$0/1$-polytopes, just as all other polytopes, can be described both
in terms 
of their vertices (``${\cal V}$-presentation'') and in terms of equations
and facet-defining inequalities (``${\cal H}$-presentation''). However,
for $0/1$-polytopes the first point of view yields the name, it gives the
natural definition, and thus it also determines our starting point. 

\begin{Definition}[0/1-polytopes]%
A {\em $0/1$-polytope\/} is a set $P\subseteq\R^d$ of the form 
\[
P=P(V)\ \ :=\ \ \conv(V)\ \ =\ \ \{V\xx:\xx\ge \zero ,\, \one^t\xx=1\}
\]
where $V\in\{0,1\}^{d\times n}$ is a $0/1$-matrix
whose set of columns, a subset of the 
vertex set of the unit cube $C_d=[0,1]^d$, is the
\emph{vertex set} of $P$.
\end{Definition}

\begin{Notation}\rm
Here and in the following, we will extensively rely on 
vector and matrix notation.
Our basic objects are column vectors such as $\xx,\yy,\ldots$. 
Their transposed vectors $\xx^t,\yy^t$ are
thus row vectors. We use $\one$ to denote a column vector of all $1$s
(whose length is defined by the context), 
$\zero$ to denote the corresponding zero vector, while $\ee_i$ denotes the
$i$-th unit vector (of unspecified length). 
The product $\xx^t\yy$ of a row with a column vector yields the
standard scalar product, while $\xx\yy^t$ is a product of a
column vector with a row vector (of the same length), and thus
represents a matrix of rank~$1$.
Thus $\one^t\one=n$ if $\one$ has length $n$, while $\one\one^t$ is 
a square all-$1$s
matrix. Matrices such as $V$ and their sets of columns are used
interchangeably. A unit matrix of size $n\times n$ will be denoted $I_n$.
\end{Notation}

It is hard to ``see'' what a $0/1$-polytope looks like from
looking at the matrix~$V$. We have more of a chance to ``understand''
an example by feeding it to a computer and asking for an analysis.
More specifically, we may present $P(V)$ to the {\tt polymake}
system of Gawrilow \& Joswig \cite{GawrilowJoswig,GawrilowJoswig2}
in terms of a file that contains the key word {\tt POINTS}
in its first line, and then the matrix $(\one,V^t)$ 
in the following lines --- the rows of this matrix give
homogeneous coordinates for the vertices of~$P(V)$.

\begin{Example}\rm
For $n\ge1$,
\begin{eqnarray*}
\Delta_{n-1} \ := \ P(I_n) &=& \conv(\{\ee_1,\ldots,\ee_n\})\\
             & = & \{\xx\in\R^n:\xx\ge\zero,\,\one^t\xx=1\}\ \ \sse\ \R^n
\end{eqnarray*}
is the \emph{standard} simplex of dimension $n-1$.\\
This is a {\em regular\/} simplex, since all its edges
have the same length~$\sqrt2$, but
it is not full-dimensional, since it lies in the hyperplane given
by $\one^t\xx=1$. 
Alternatively, we could consider the simplex
\begin{eqnarray*}
\Delta'_n \ = \ P(\zero,I_n) &=& \conv(\{\zero,\ee_1,\ldots,\ee_n\})\\
             & = & \{\xx\in\R^n:\xx\ge\zero,\,\one^t\xx\le1\}\ \ \sse\ \R^n,
\end{eqnarray*}
which is full-dimensional, but not regular for $n\ge2$.
In fact, in many dimensions (starting at~$n=2$)
there is no full-dimensional, regular $0/1$-simplex at all.
(See Problem~\ref{prob:reg_sx}.)
\end{Example}

\begin{ExampleExercise}
For $V\in\{0,1\}^{d\times n}$, let 
$\widetilde V=\binom{V_{~}}{I_n}\in\{0,1\}^{(d+n)\times n}$. Then
\[
\begin{array}{ll}
P(\widetilde V)\ =\ \Big\{\,\dbinom{\xx}{\yy}\in[0,1]^{d+n}: 
& \yy\ge 0,\ \one^t\yy=1,\\
                  & x_i=\sum\limits^n_{j=1} v_{ij}y_j\textrm{ for
}1\le
i\le d\,\Big\}
\end{array}
\]
is an affine image of the $(n-1)$-dimensional standard simplex~$\Delta_{n-1}$.
(Prove this!)
Thus for the $0/1$-polytope $P(\widetilde V)\sse\R^{d+n}$ we
\emph{have} a complete description in terms of linear equations and
inequalities. From this we get $P(V)$ as the image of the projection
\begin{eqnarray*}
\pi: \ \R^{d+n} & \longrightarrow & \R^d\\
        \dbinom{\xx}{\yy} & \longmapsto &  \xx
\end{eqnarray*}
that deletes the last $n$ coordinates. Equivalently, 
to get $P(V)=\pi(P(\widetilde V))$ from $P(\widetilde V)$
we must apply the
operation ``delete the last coordinate'' $n$ times. 
\end{ExampleExercise}

\begin{Theorem}[$\cal H$-presentations]%
Every $0/1$-polytope $P(V)\sse\R^d$ 
can be written as the set of solutions of a
system of linear inequalities, that is, as
\[
P(V)=\{\xx\in\R^d:A\xx\le \bb\}
\]
for some $n\in\N$, a matrix $A\in\Z^{n\times d}$, and a vector
$\bb\in\Z^n$.
\end{Theorem}

\proof
First, we need not deal with 
equations in the system that describes $P(\widetilde V)$, 
since these can be rewritten in terms of
inequalities: the equation $\aaa^t\xx= \beta$ is equivalent to
the two inequalities $\aaa^t\xx\le \beta$, $-\aaa^t\xx\le -\beta$.
Thus, with the observations above, it suffices to show that if a set
$S\subseteq \R^{k+1}$ has a description of the form 
\[
S\ \ =\ \ \Big\{\dbinom{\xx}{x_{k+1}}\in\R^{k+1}:
\aaa^t_i\xx+a_{i,k+1}x_{k+1}\le b_i\ \ (1\le i\le m)\Big\},
\]
then the projection of $S$ to $\pi(S)\subseteq \R^k$ (by ``deleting
the last coordinate'') has a representation of the same type. We may
assume that the inequality system has been ordered so that 
\begin{eqnarray*}%
a_{i,k+1} > 0& \textrm{ for } &1\ \le i\le i_0,\\
a_{j,k+1} < 0& \textrm{ for } &i_0<j\le j_0,\\
a_{i,k+1} = 0& \textrm{ for } &j_0<i\le m.
\end{eqnarray*}
Now for any given $\xx\in\R^k$, it is easy to decide whether
it lies in $\pi(S)$. Namely, $\xx\in\pi(S) $ holds if and only if
there is some value $\xi\in\R$ such that $\binom\smallxx\xi\in S$,
where the inequalities for $1\le i\le i_0$ provide upper bounds 
for such a value~$\xi$, the inequalities for $i_0<j\le j_0$
give lower bounds, the others provide no conditions.
Thus the system has a solution $\xi$ for given $\xx$ if
all the upper bounds are at least as large as all the lower bounds. 
Explicitly, this yields a description of $\pi(S)$ as
\[
\left\{
\begin{array}{clcll}
\xx\in\R^k:&(a_{i,k+1}\aaa_j
-a_{j,k+1}\aaa_i)^t\xx&\le& a_{i,k+1}b_j-a_{j,k+1}b_i
&\textrm{for }1\le i\le i_0
\\ &&&&\textrm{ \ \ and }\ i_0<j\le j_0,\\
&a_i\xx&\le& b_i
&\textrm{for }j_0<i\le m
\end{array}\right\},
\]
which is a presentation of the required form.
\proofend

The transformation of an inequality system for $S$ into a system for
$\pi(S)$ in this way is known as \emph{Fourier-Motzkin elimination}
of the last variable \cite[Lecture~1]{Zpoly}.
Note that, in the worst case, the system for $\pi(S)$
may have as many as $\left(\frac{m}{2}\right)^2$ inequalities:
much more than the system for $S$!
The good news at this point is that the inequality descriptions
of $\pi(S)$ are typically very redundant: many of the inequalities
can be deleted without changing the set of solutions of the system.
However, the bad news is that even a minimal system --- which
in the case of a full-dimensional polytope $P$ consists
of exactly one inequality for each facet of~$P$ --- may be huge.
Correspondingly, $0/1$-polytopes with rather few vertices
may have ``many'' facets: See Section~\ref{sec:facets} below.

A projection argument together with the basic operation of
``switching'' will allow us for the following
to assume that the polytopes under consideration are
full-dimensional, and have $\zero$ as a vertex, whenever that seems
convenient:
\begin{itemize}
\item[(1)] All the symmetries of the $0/1$-cube $C_d=[0,1]^d$ transform
$0/1$-polytopes into $0/1$-polytopes. In coordinates, these 
symmetries are generated by\\
$\bullet$ permuting coordinates, and\\
$\bullet$ replacing some coordinates $x_i$ by $\overline{x}_i:=1-x_i$
({\it switching\/}).
\[
\begin{picture}(0,0)%
\epsfig{file=EPS/switch.pstex}%
\end{picture}%
\setlength{\unitlength}{2763sp}%
\begingroup\makeatletter\ifx\SetFigFont\undefined%
\gdef\SetFigFont#1#2#3#4#5{%
  \reset@font\fontsize{#1}{#2pt}%
  \fontfamily{#3}\fontseries{#4}\fontshape{#5}%
  \selectfont}%
\fi\endgroup%
\begin{picture}(6040,1391)(1118,-4558)
\put(2626,-3661){\makebox(0,0)[lb]{\smash{\SetFigFont{12}{14.4}{\rmdefault}{\mddefault}{\updefault}$\sim$}}}
\put(5551,-4486){\makebox(0,0)[b]{\smash{\SetFigFont{12}{14.4}{\rmdefault}{\mddefault}{\updefault}$x_1\longleftrightarrow1-x_1$}}}
\put(2626,-4486){\makebox(0,0)[b]{\smash{\SetFigFont{12}{14.4}{\rmdefault}{\mddefault}{\updefault}$x_1\longleftrightarrow x_2$}}}
\put(5176,-3661){\makebox(0,0)[lb]{\smash{\SetFigFont{12}{14.4}{\rmdefault}{\mddefault}{\updefault}$\sim$}}}
\end{picture}

\]
We call two $0/1$-polytopes $P$ and $P'$ \emph{$0/1$-equivalent} if 
a sequence of such operations
can transform $P$ into $P'$. In particular, one can 
transform any $0/1$-polytope $P$ with a vertex 
$\vv\in P\cap\{0,1\}^n$ to a new, $0/1$-equivalent polytope $P'$
such that the vertex $\vv$ gets mapped to the vertex $\zero $ of~$P'$.
\item[(2)] If $P\subseteq \R^{d+1}$ is not \emph{full-dimensional},
then it is affinely equivalent to a $0/1$-polytope $P'\subseteq\R^d$. 
To see this, first we may assume that $\zero \in P$ (after switching),
so $P$ satisfies an equation of the form $\aaa^t\xx+a_{d+1}x_{d+1}=0$ 
with $\aaa\in\R^d $. 
\[
\begin{picture}(0,0)%
\epsfig{file=EPS/proj.pstex}%
\end{picture}%
\setlength{\unitlength}{2763sp}%
\begingroup\makeatletter\ifx\SetFigFont\undefined%
\gdef\SetFigFont#1#2#3#4#5{%
  \reset@font\fontsize{#1}{#2pt}%
  \fontfamily{#3}\fontseries{#4}\fontshape{#5}%
  \selectfont}%
\fi\endgroup%
\begin{picture}(2707,1469)(5172,-4658)
\put(6526,-4336){\makebox(0,0)[b]{\smash{\SetFigFont{12}{14.4}{\rmdefault}{\mddefault}{\updefault}$\downarrow\downarrow$}}}
\end{picture}

\]
Furthermore, after
permuting the coordinates we get that $a_{d+1}\not= 0$. But then
``deleting the last coordinate''
\[
\pi:\R^{d+1}\rightarrow\R^d
\]
projects $P\rightarrow P'=\pi(P)$ injectively, that is, it defines an
affine equivalence between $P$ and $\pi(P)=P'$.
\end{itemize}

In the following, we usually deal with full-dimensional
$0/1$-polytopes, and we take $0/1$-equivalence as the basic
notion for their comparison. 
The resulting classification is much finer than the classification
by affine equivalence --- for example, all $d$-dimensional 
$0/1$-simplices are affinely equivalent, but they are not necessarily
$0/1$-equivalent: Note that
$0/1$-equivalent polytopes are congruent, so they have the
same edge lengths, volumes, etc. But the converse is not true,
see below.

\begin{Definition}%
The \emph{faces} of a $0/1$-polytope  $P$ are the subsets of the
form $P\supcc=\{\xx\in P:\cc^t\xx=\gamma\}$, where $\cc^t\xx\le\gamma$
is a linear inequality that is valid for \emph{all} points of~$P$.
This definition of faces includes the subsets $\emptyset$ and $P$, the 
\emph{trivial faces} of~$P$.

All faces of a $0/1$-polytope are themselves $0/1$-polytopes, of the
form $F=\conv(F\cap\{0,1\}^d)$.
The set of $0$-dimensional faces, or \emph{vertices}, of a $0/1$-polytope
is given by $V=P\cap\Z^d$. The $1$-dimensional faces are called \emph{edges}.
Vertices and edges together form the \emph{graph} of the polytope.   
The maximal non-trivial faces, of dimension $\dim(P)-1$,
are the \emph{facets} of~$P$.
These are essential for the $\cal H$-presentation of polytopes:
In the full-dimensional case every irredundant $\cal H$-presentation
consists of exactly one inequality for each facet of~$P$. 

The \emph{face lattice} is the set of all faces of~$P$, 
partially ordered by inclusion. 
It is a graded lattice of length $\dim(P)+1$.
Two polytopes are \emph{combinatorially equivalent} 
if their face lattices are isomorphic as finite lattices. 
\end{Definition}

\begin{Proposition}%
On the finite set of all $0/1$-polytopes in~$\R^d$
one has the following hierarchy of equivalence relations:
\begin{center}
``$0/1$-equivalent''$\ \Rightarrow\ $``congruent''$\ \Rightarrow\ $``affinely
equivalent''$\ \Rightarrow\ $``combinatorially equivalent.''
\end{center}
For all three implications the converse is false, even when we restrict
the discussion to full-dimensional polytopes.
\end{Proposition}

\proof
The hierarchy is clearly valid: Every $0/1$-equivalence is a congruence,
congruent polyhedra are affinely equivalent, and affine equivalence
implies combinatorial equivalence. In the following we provide
counterexamples for all the converse implications.

(1)
Full-dimensional $0/1$-polytopes that are
congruent but not $0/1$-equivalent can be found in
dimension~$5$:
\begin{verbatim}
VERTICES              VERTICES
1 0 0 0 0 0           1 0 0 0 0 0
1 0 0 1 1 0           1 0 0 1 1 0
1 0 1 0 1 0           1 0 1 0 1 0
1 1 0 0 1 0           1 0 1 1 0 0
1 0 1 1 0 0           1 1 0 0 1 0
1 0 1 1 0 1           1 1 0 0 1 1
\end{verbatim}
One easily checks that these two data sets (in {\tt polymake} input format;
see {\tt CNG:5-6a.poly} 
and {\tt CNG:5-6b.poly} in the {\tt polymake} database)
describe congruent, full-dimensional $0/1$-simplices in~$\R^5$:
For this one just computes the pairwise distances of the points.
A $0/1$-equivalence would transform the array on the left to the
array on the right by permuting rows and columns, and by complementing
columns. But on the left we have two columns with exactly one~$1$
(and no column with five $1$s), while on the right there
is only one column with exactly one~$1$ (and no column with five $1$s).
Thus the two simplices are not $0/1$-equivalent.
Volker Kaibel has additionally shown that for $d\le4$ all congruent
full-dimensional 
$0/1$-polytopes are indeed $0/1$-equivalent.

(2)
The above classification for $d=3$ contains examples of tetrahedra 
that are not congruent, but of course affinely equivalent.
Further examples will appear in Lecture~2.

(3)
Here are two $5$-polytopes,
{\tt EQU:5-7a.poly} and {\tt EQU:5-7b.poly}, that are
combinatorially, but not affinely equivalent:
\begin{verbatim}
VERTICES              VERTICES
1 0 0 0 0 0           1 0 0 0 0 0
1 1 0 0 0 0           1 1 1 0 0 0
1 0 1 0 0 0           1 0 1 1 0 0
1 0 0 1 0 0           1 0 0 1 1 0
1 0 0 0 1 0           1 0 0 0 1 1
1 0 0 0 0 1           1 1 0 0 0 1
1 1 1 1 1 1           1 1 1 1 1 1
\end{verbatim}
In fact, each of them is a bipyramid
over a 4-simplex (and hence they are combinatorially
equivalent), but in the first one the main diagonal
is divided in the ratio $1:4$, 
for the other one the ratio is $2:3$,
and such ratios are preserved by affine equivalences.
\proofend


\subsection{Doubly exponentially many $0/1$-polytopes}

How many non-equivalent $0/1$-polytopes are there?
Clearly in $\R^d$ there are exactly $2^{2^d}$
different $0/1$-polytopes, but some of them are
low-dimensional, and some of them are equivalent to many others.
Nevertheless, this trivial estimate is not that far from 
the truth.

For the following, let $F_i^0$ denote the facet of the
$d$-cube $[0,1]^d$ that is given by $x_i=0$, and similarly let
$F_i^1$ be the facet given by $x_i=1$.
With a $3$-dimensional picture in the back of our minds,
we will refer to $F_d^0$ as the {\em bottom facet\/} and to
$F_d^1$ as the {\em top facet\/} of~$C_d$.
All other facets will be called the {\em vertical facets\/}
of~$C_d$.

This terminology corresponds to one of the main proof 
techniques that we have for $0/1$-polytopes: decomposition into
``top'' and ``bottom'' with induction over the dimension.
For this we note the following for an arbitrary
$0/1$-polytope $P\sse[0,1]^d$:
\begin{itemize}
\item Every facet $F_i^s$ induces a face $P_i^s:=F_i^s\cap P$ of~$P$;
         these faces are referred to as the {\em trivial faces} of~$P$.
\item Every vertex of~$P$ is contained either in the 
   {\em bottom face\/} $P_d^0=F_d^0\cap P$ or in the
   {\em top face\/} $P_d^1=F_d^1\cap P$ of~$P$.
\item Every vertex $\vv$ of~$P$ is determined by the set of trivial
   faces~$P_i^0$ that contain it, since
   $v_i=0$ holds if and only if $\vv\in P_i^0$. 
\end{itemize}
The following figure illustrates that in general some
trivial faces are facets, while others are not.
\[
\begin{picture}(0,0)%
\epsfig{file=EPS/faces.pstex}%
\end{picture}%
\setlength{\unitlength}{3947sp}%
\begingroup\makeatletter\ifx\SetFigFont\undefined%
\gdef\SetFigFont#1#2#3#4#5{%
  \reset@font\fontsize{#1}{#2pt}%
  \fontfamily{#3}\fontseries{#4}\fontshape{#5}%
  \selectfont}%
\fi\endgroup%
\begin{picture}(3512,2331)(4439,-7111)
\put(6226,-4936){\makebox(0,0)[rb]{\smash{\SetFigFont{12}{14.4}{\rmdefault}{\mddefault}{\updefault}top face $P_3^1$}}}
\put(6826,-7111){\makebox(0,0)[lb]{\smash{\SetFigFont{12}{14.4}{\rmdefault}{\mddefault}{\updefault}bottom face $P_3^0$}}}
\put(6826,-6511){\makebox(0,0)[lb]{\smash{\SetFigFont{12}{14.4}{\rmdefault}{\mddefault}{\updefault}$P$}}}
\put(7951,-5086){\makebox(0,0)[lb]{\smash{\SetFigFont{12}{14.4}{\rmdefault}{\mddefault}{\updefault}$P_2^1$}}}
\put(5176,-5761){\makebox(0,0)[rb]{\smash{\SetFigFont{12}{14.4}{\rmdefault}{\mddefault}{\updefault}$P_1^0$}}}
\put(5101,-7111){\makebox(0,0)[rb]{\smash{\SetFigFont{12}{14.4}{\rmdefault}{\mddefault}{\updefault}$P_2^0$}}}
\put(7876,-6061){\makebox(0,0)[lb]{\smash{\SetFigFont{12}{14.4}{\rmdefault}{\mddefault}{\updefault}$P_1^1$}}}
\end{picture}

\]
\smallskip

\begin{Proposition}[Sarangarajan-Ziegler]
There is a family ${\cal F}_d$ of $2^{2^{d-1}-4}$ different,
full-dimensional $0/1$-polytopes in~$[0,1]^d$, such that
\begin{itemize}
\item
any two polytopes in ${\cal F}_d$ are $0/1$-equivalent 
if and only if they are combinatorially equivalent, and
\item
for $d\ge6$, the collection ${\cal F}_d$ contains 
more than $2^{2^{d-2}}$ combinatorially non-equivalent
$d$-dimensional $0/1$-polytopes in~$\R^d$.
\end{itemize}
\end{Proposition}

\proof
Let $d\ge3$, and let ${\cal F}_d$ be
the set of $0/1$-polytopes $P(V)=\conv(V)$ 
of the following form:
\begin{itemize}
\item
$V$ contains all the vertices in the bottom facet $F_d^0$ of the
$d$-cube $[0,1]^d$ (that is, $\{0,1\}^{d-1}\times\{0\}\sse V$), 
\item
the pair $\ee_d,\one$ of opposite vertices of the top facet $F_d^1$
is contained in~$V$,
\item
the two opposite vertices $\ee_d+\ee_1,\one-\ee_1$ of the top facet $F_d^1$
are {\em not\/} contained in~$V$.
\end{itemize}
This fixes $2^{d-1}+4$ vertices to be inside or outside~$V$, and thus
leaves $2^{2^d-(2^{d-1}+4)}=2^{2^{d-1}-4}$ 
choices for the set~$V$, and hence for the polytope~$P(V)$.

For $d=3$, there is exactly one polytope of the
given special type (the ``nameless'' one):
\[
\epsffile{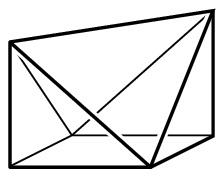}
\]
Now the following facts are easy to verify about the
polytopes~$P(V)\in{\cal F}_d$:
\begin{itemize}
\item 
$P(V)$ is a $d$-dimensional $0/1$-polytope.
Its bottom facet $P_d^0=F_d^0$ is a $(d-1)$-cube,
with $2^{d-1}$ vertices.
\item
All the vertical facets $F_i^s$ ($i<d$) induce facets $P_i^s$
of~$P(V)$. 
These are the facets of~$P(V)$ that are adjacent to the cube facet
$P_d^0$.
Every vertex of $P(V)$ that is not on $F_d^0$ is completely determined by
the set of vertical facets $P_i^s$ that it lies on.
\item
All facets of $P(V)$, other than the bottom
facet, have fewer than~$2^{d-1}$ vertices.
\\
(For this we use that only the $2d+2\binom{d}{2}=d^2+d$
``special'' hyperplanes given by 
$x_i=0$, $x_i=1$ or $x_i=x_j$ or $x_i=1-x_j$ contain $2^{d-1}$ 
$0/1$-points, and all other hyperplanes
contain less than $2^{d-1}$ $0/1$-points. 
It is easy to verify that no special hyperplane other than ``$x_d=0$''
can describe a facet of $P(V)$.)
\item
Therefore, if two polytopes $P(V)$ and $P(V')$ are
combinatorially isomorphic, then they are equivalent by a
symmetry of the $d$-dimensional $0/1$-cube that fixes the
bottom facet, and induces an automorphism of that bottom facet.
\item
The order of the symmetry group of $C_{d-1}$ is $2^{d-1}(d{-}1)!$.
So for each $P(V)$ there are not more than $ 2^{d-1}(d{-}1)!$
polytopes $P(V')$ that are combinatorially equivalent to it. 
\item
Therefore, there are more than
$2^{2^{d-1}-4}/(2^{d-1}(d{-}1)!)$ combinatorially non-isomorphic 
$0/1$-polytopes of the form $P(V)$, and for
$d>5$ this number is larger than $2^{2^{d-2}}$. 
\proofend
\end{itemize}

\section{The Number of Facets}\label{sec:facets}

\subsection{Some examples}

Staring too much at the $3$-dimensional case, one might come
up with the conjecture that a $d$-dimensional $0/1$-polytope
cannot have more than $2^d$ facets.
In fact, 
\[
C_d^\Delta\ \ :=\ \ \conv\{\ee_1,\ldots,\ee_d,\one-\ee_1,\ldots,\one-\ee_d\}
\]
is a polytope with $2d$ vertices ($d\ge3$) that is
centrally symmetric with respect to $\frac12\one$,
the center of the $0/1$-cube. Hence it is
affinely equivalent to the usual regular $d$-dimensional cross
polytope. In particular, this polytope has $2^d$ facets.
The first examples are given as
{\tt CRO:3-6.poly}, {\tt CRO:4-8.poly}, \ldots
in the database.

(For $d=4$ this construction produces a regular cross
polytope {\tt CRO:4-8.poly}, all of whose edges have length $\sqrt2$.
Another remarkable regular cross polytope {\tt HAM:8-16.poly}
arises from the extended Hamming code $\widetilde{H}_8$.
The cross polytopes $C_d^\Delta$ as constructed above are not
regular for $d\neq4$: they have edges of lengths $\sqrt2$ and $\sqrt{d-2}$.)
\newcommand{\facet}{\operatorname{\#f}}

But more than that? Ewgenij Gawrilow was the first to detect a 
$5$-dimensional $0/1$-polytope with $40$ facets. 
After intensive search, here is what we know about examples
of low-dimensional $0/1$-polytopes with ``many facets'' ---
and thus about $\facet(d)$,
the maximal number $f_{d-1}(P)$ of facets that a $d$-dimensional
$0/1$-polytope $P$ can have:
\[
\begin{array}{c|r@{\quad}lr}
d & \facet(d)\ \qquad  & \mbox{proved/found by}\!\! & \mbox{example}\\
\hline \\[-4mm]
3 & ~= 8 &&  \mbox{\tt CRO:3-6.poly}\\
4 & ~ = 16 & \mbox{Below} & \mbox{\tt CRO:4-8.poly}\\
5 & ~ = 40 & \mbox{Aichholzer} & \mbox{\tt EG:5-10.poly}\\
6 & ~ \ge 121 & \mbox{Sarangarajan} & \mbox{\tt AS:6-18.poly}\\
7 & ~ \ge 432 & \mbox{Christof} & \mbox{\tt TC:7-30.poly}\\
8 & ~ \ge 1675 & \mbox{Christof} & \mbox{\tt TC:8-38.poly}\\
9 & ~ \ge 6875 & \mbox{Christof} & \mbox{\tt TC:9-48.poly}\\
10 & ~ \ge 41591 & \mbox{Christof} & \mbox{\tt TC:10-83.poly} \\
\vdots && \quad\quad\vdots & \vdots\quad\quad \\
13 & ~ \ge  17464356 & \mbox{Christof} & \mbox{\tt TC:13-254.poly}
\end{array}
\]
In brief: $0/1$-polytopes may have {\em many\/} facets.
But how many, at most? And how do $0/1$-polytopes with 
``many facets'' look like?


\subsection{Some upper bounds}

The asymptotically best upper bound for the number 
of facets of a $d$-dimensional $0/1$-polytope is 
the following. I assume that it is rather tight;
the problem is with the lower bounds, which look 
much worse.

\begin{Theorem}[Fleiner, Kaibel \& Rote \cite{FleinerKaibelRote}]%
For all large enough $d$, a $d$-dimensional $0/1$-polytope
has no more than 
\[
\facet(d)\ \le\ 30\,(d-2)!
\]
facets.
\end{Theorem}

See \cite{FleinerKaibelRote} for the (beautiful) proof of this result,
which is probably valid for {\em all\/}~$d$.
The first bound of this order of magnitude was pointed out by
Imre B\'ar\'any \cite[p.~26]{Zpoly}.
Here we present a proof for the inequality
\begin{equation}\tag{$*$}
\facet(d)\ \le\ 2\,(d-1)! \ +\ 2(d-1),
\end{equation}
which is asymptotically a bit worse than the one just quoted, but
it is better in low dimensions --- and whose proof 
(also from \cite{FleinerKaibelRote}) is strikingly simple.

For this, let $P\sse[0,1]^d$ be a $d$-dimensional $0/1$-polytope.
We note the following facts:
\begin{itemize}
\item The volume $\Vol_d(P)$ is an integral multiple of $\frac1{d!}$.\\
    (Every polytope can be triangulated without new vertices.
    Thus we are reduced to the case of $0/1$-simplices, whose
    volume is given as $\frac1{d!}$ times the determinant --- which
    is an integer.)
\item \label{eq:vol}
    The number of facets $f_{d-1}(P)$ 
    of a $d$-dimensional $0/1$-polytope $P$ satisfies 
\[ 
f_{d-1}(P)\ \ \le\ \  2d \ +\ d!\,(1-\Vol_d(P)).
\] 
    (This follows from an observation of B\'ar\'any:
    The $d$-cube $[0,1]^d$ has $2d$ facets. Now delete the ``superfluous''
    $0/1$-vectors, so that $[0,1]^d$ is gradually transformed into $P$. 
    Whenever a facet of $P$ ``appears'' in this process,
    a pyramid over the facet is removed,
    and the volume of this pyramid is at least $\frac1{d!}$.)
\item
    Consider the projection $\pi:\R^d\longrightarrow\R^{d-1}$ 
    that deletes the last 
    coordinate. With respect to this projection, the boundary of~$P$ may be
    divided into ```vertical,'' ``upper'' and ``lower'' facets.
    After projection, the images of the upper facets partition 
    $\pi(P)$ into $(d-1)$-dimensional $0/1$-polytopes, and so
    do the lower facets. Thus we get that 
\newcommand{\upper}{\operatorname{upper}}
\newcommand{\low}{\operatorname{lower}}
\newcommand{\ver}{\operatorname{vertical}}
\[
f^{\upper}_{d-1}(P),\,
f^{\low}_{d-1}(P) \ \ \le \ \ (d-1)!\, \Vol_{d-1}(\pi(P)).
\]
Our figure illustrates this decomposition for $d=3$:
\[
\begin{picture}(0,0)%
\epsfig{file=EPS/vert.pstex}%
\end{picture}%
\setlength{\unitlength}{3947sp}%
\begingroup\makeatletter\ifx\SetFigFont\undefined%
\gdef\SetFigFont#1#2#3#4#5{%
  \reset@font\fontsize{#1}{#2pt}%
  \fontfamily{#3}\fontseries{#4}\fontshape{#5}%
  \selectfont}%
\fi\endgroup%
\begin{picture}(2858,1694)(4343,-6833)
\put(7201,-6811){\makebox(0,0)[lb]{\smash{\SetFigFont{12}{14.4}{\rmdefault}{\mddefault}{\updefault}$f^{\low}_2=2$}}}
\put(7201,-5311){\makebox(0,0)[lb]{\smash{\SetFigFont{12}{14.4}{\rmdefault}{\mddefault}{\updefault}$f^{\upper}_2=2$}}}
\put(7201,-6136){\makebox(0,0)[lb]{\smash{\SetFigFont{12}{14.4}{\rmdefault}{\mddefault}{\updefault}$P$}}}
\put(4951,-5686){\makebox(0,0)[b]{\smash{\SetFigFont{12}{14.4}{\rmdefault}{\mddefault}{\updefault}$f^{\ver}_2=2$}}}
\end{picture}

\]
\item
    At the same time, the vertical facets of~$P$ are in bijection
    with a subset of the facets of~$\pi(P)$: and the number of these can be
    estimated using the formula above:
\[
f^{\ver}_{d-1}(P)\ \ \le\ \ f_{d-2}(\pi(P))
                 \ \ \le\ \ 2(d-1) + (d-1)!\,(1-\Vol_{d-1}(\pi(P))).
\]
\item
    \ldots and summing the upper bounds that we have obtained
    for $f^{\upper}_{d-1}(P)$, $f^{\low}_{d-1}(P)$ and
    $f^{\ver}_{d-1}(P)$ completes the proof of~$(*)$.
\proofend
\end{itemize}


\subsection{A bad construction}

All the available data suggests that 
$0/1$-polytopes may have much more than just simply-exponentially
many facets.
But no one has been able, up to now, to prove any lower bound
on $\facet(d)$ that grows faster than $c^d$, for some constant $c>1$.

\begin{Proposition}[Kortenkamp et al.~\cite{KRSZ}]%
For all large enough $d$,
\[
\facet(d)\ >\  3.6^d.
\]
\end{Proposition}

\proof
The {\em sum\/} $P_1*P_2$  of two polytopes $P_1$ and $P_2$
is obtained by representing the polytopes in some $\R^n$ in such a way 
that their intersection
consists of one single point which
for both of them lies in the relative interior,
and by then taking the convex hull:
\[
P\ :=\ \conv(P_1\cup P_2),\quad\textrm{if $P_1\cap P_2=\{\xx\}$ is
a relative interior point for both  $P_1$ and $P_2$.}
\]
If we take the sum of two polytopes in this way, then 
the dimensions add, while the
number of facets are multiplied.
As an example, the sum of an $n$-gon (dimension~$2$, $n$~facets)
and an interval (dimension~$1$, $2$~facets)
results in a bipyramid over the $n$-gon (dimension~$3$, $2n$~facets).
The sum operation is polar to taking products, where the
dimensions add and the numbers of vertices are multiplied.

But we have to take a bit of care in order to 
adapt this general polytope operation to $0/1$-polytopes,
since there is very little space for ``moving into a position'' if
we want to stay within the setting of $0/1$-polytopes.
For this call a $0/1$-polytope {\em centered\/} if it has
the center point $\frac12\one$ in its (relative) interior.
For example, among the $3$-dimensional $0/1$-polytopes,
the $3$-dimensional prism, the two different
pyramids over a square, and the tetrahedra except for $\CUT(3)$, 
are {\em not\/} centered! 
On the other hand, the cross polytopes $C_d^\Delta$ are centered
for all~$d$.

The sum of two centered $0/1$-polytopes $P_1\sse[0,1]^{d_1}$ and
$P_2\sse[0,1]^{d_2}$ can be realized 
in $[0,1]^{d_1+d_2}$ by embedding them into the subspaces 
$x_{d_1}=x_{d_1+1}=\ldots=x_{d_1+d_2}$
resp.\
$x_1=\ldots=x_{d_1}=1-x_{d_1+1}$.
This yields centered
$0/1$-polytopes $\widehat{P}_1,\widehat{P}_2\sse[0,1]^{d_1+d_2}$
that are affinely isomorphic to $P_1$ and $P_2$, and whose
convex hull realizes the sum $P_1*P_2$.
For example, the octahedron $C_3^\Delta$
can be viewed as the sum of a rectangle
and a diagonal:
\[
\epsffile{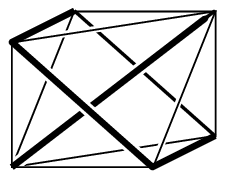}
\]

Now we need a starting block: and for that we use
Christof's $13$-dimensional $0/1$-polytope {\tt TC:13-254.poly}
with at least $17464356 > 3.6^{13}$ facets.
This polytope is indeed centered (you may check that already the
first $22$ vertices contain $\frac12\one$ in their interior).
Taking sums of copies of this polytope, 
and extra copies of $[0,1]$ if needed,
we arrive at the result.
\proofend

This seems to be the best asymptotic lower bound available in the moment.
I think that it is {\em bad\/}: one should be able to prove a lower bound
of the form $c^{d\log d}$, or at least
that there is a lower bound that grows faster than $C^d$ for
every $C>1$!
I'd offer two candidates for such a lower bound construction:
Random polytopes, and cut polytopes. However, we cannot do the
corresponding lower bound proof for either of these two classes,
up to now.


\section{Random $0/1$-Polytopes}

We do not understand random $0/1$-polytopes very well.
Let $d$ be not too small, and take, say, $2d$ or $d\log d$ or $d^2$ random 
$0/1$-points:
{\em How will their convex hull look like? How many edges, and how
many facets can we expect the random polytope to have?}
We will see in this section that the analysis of 
random $0/1$-polytopes is driven by one basic linear algebra parameter:
the probability $P_d$ that a random $0/1$-matrix of size $d\times d$ 
has vanishing determinant.

This probability corresponds to the case of $d+1$ random $0/1$-points:
Take $d+1$ points $\vv_0,\vv_1,\ldots,\vv_d$ independently
at random (where all $0/1$ points appear with the same probability
$p=\frac1{2^d}$).
The $d+1$ points will be distinct with very high probability,
and by symmetry we may assume that the first point is $\vv_0=\zero$.
Thus the probability that the $d+1$ points span a $d$-dimensional
simplex is exactly~$1-P_d$. How large is this probability?
We first study the case where $d$ is small, and from this
we will derive a quite misleading impression.


\subsection{The determinant of a small random $0/1$-matrix}

Let $P_d$ be the probability that a random $0/1$-matrix
of size~$d\times d$ is singular. Of course we have
\[
P_d\ =\frac{M_d}{2^{d^2}},
\]
where $M_d$ denotes the number of different $0/1$-matrices
of size $d\times d$ that have determinant~$0$. 
This number can be computed explicitly for $d\le7$:
\[
\begin{array}{r|r}
d & M_d\\
\hline
\\[-4mm]
1 &               1 \\
2 &              10 \\
3 &             338 \\
4 &           42976 \\
5 &        21040112 \\
6 &     39882864736 \\
7 & 292604283435872
\end{array}
\]
(In fact, for $d\le6$ numbers that are equivalent to these
were computed by Mark B. Wells in the sixties \cite[p.~198]{MetropolisStein};
the value for $d=7$ is new, due to Gerald Stein.)

From this, we get a table for $P_d$, where for $d\ge8$ we print
estimates that were obtained by taking 10 million random matrices for each
case:

\[
\begin{array}{lcl}
P_1 & = & 0.5 \\[-.3mm]
P_2 & = & 0.625 \\[-.3mm]
P_3 & = & 0.66015625 \\[-.3mm]
P_4 & = & 0.65576... \\[-.3mm]
P_5 & = & 0.62704... \\[-.3mm]
P_6 & = & 0.58037... \\[-.3mm]
P_7 & = & 0.51976... \\[-.3mm]
P_8 &\approx & 0.449 \\[-.3mm]
P_9 &\approx & 0.373 \\[-.3mm]
P_{10} &\approx & 0.298\\[-.3mm]
P_{11} &\approx & 0.226\\[-.3mm]
P_{12} &\approx & 0.164\\[-.3mm]
P_{13} &\approx & 0.113\\[-.3mm]
P_{14} &\approx & 0.075\\[-.3mm]
P_{15} &\approx & 0.047 
\end{array}
\]
Conclusion: the probability $P_d$ first increases (!), but then
it seems to decrease and approach~$0$ steadily, but 
not very fast.




\subsection{Koml\'os' theorem and its consequences}

The question about the probability $P_d$
of singular random $0/1$-matrices is
equivalent to the same question about $\pm1$-matrices:
$P_d$ is equally the probability that a random
$\pm1$-matrix of size $(d+1)\times(d+1)$ is singular.
It is often convenient to switch to 
the $\pm1$-case because it has more symmetry.
The following proposition 
establishes the equivalence. Its observation is quite trivial, but also
fundamental for various
problems related to $0/1$-polytopes.

\begin{Proposition}[Williamson~\cite{Williamson}]\label{Prop:biject}%
The map 
\[
\varphi:\ A\ \ \longmapsto\ \ 
\left(
\begin{array}{cc}
1       & \one^t\\
\one & \one\one^t-2A
\end{array}
\right)
\ \ =:\ \ \Ahat.
\]
establishes a bijection between the
$0/1$-matrices of size $d\times d$ and the
$\pm1$-matrices of size $(d+1)\times(d+1)$
for which all entries in the first row and column are~$+1$.

The bijection $\varphi$ satisfies $\det(\Ahat) = (-2)^d\,\det(A)$.
In particular, it also provides a bijection between
the invertible matrices of the two types.

Furthermore, there is a one-to-$2^{2d+1}$ correspondence
between the $0/1$-matrices of size $d\times d$ and the
$\pm1$-matrices of size $(d+1)\times(d+1)$.
The correspondence again respects invertibility.
\end{Proposition}

\proof
Geometrically, the map $\varphi$
realizes an embedding of $[0,1]^d$ as a facet
of $[-1,+1]^{d+1}$.
Algebraically, 
$\Ahat =
\left(
\begin{array}{cc}
1            & \zero^t\\
\one  & I_d
\end{array}
\right)
\left(
\begin{array}{cc}
1            & \one^t\\
\zero  & -2A
\end{array}
\right)$ arises from
$\left(
\begin{array}{cc}
1            & \one^t\\
\zero  & -2A
\end{array}
\right)$
by adding the first row to all others, and thus we see that
$\varphi(A)$ is indeed invertible if $A$ is, and that
\ $\det(\Ahat)=(-2)^d\det(A)$.

Finally, with every $\pm1$-matrix 
one can associate a canonical matrix of the same size and type
for which the first row and column are positive: for this
first multiply columns by $-1$ in order to make the first row
positive, then multiply rows to make the first column positive.
There are exactly $2^{2d+1}$ matrices in
$\{-1,+1\}^{(d+1)\times(d+1)}$ that have the same canonical
form, corresponding to the $2d+1$ entries in the first row and
column for which a sign can be chosen.
\proofend

Thus $P_d$ measures for $0/1$-matrices as well as for
$\pm1$-matrices the probability of determinant~$0$. Our
experimental evidence is that $P_d$ should converge to~$0$.
But how fast?
Here is what we know.

\begin{Theorem}[Koml\'os' Theorem; Kahn, Koml\'os and Szemer\'edi \cite{KKS}]%
\label{thm:komlos}%
The probability $P_d$ that a random $0/1$-matrix of size $d\times d$ is 
singular satisfies
\[
\frac{d^2}{2^d}\ \ <\ \ P_d\ \ <\ \ 0.999^d
\]
for all high enough~$d$.
\end{Theorem}

\proof
The non-trivial part is the upper bound, which
is due to Kahn, Koml\'os and Szemer\'edi \cite{KKS}.
Their proof is difficult, involving a probabilistic construction.
In fact, it is hard enough to prove that $P_d$ converges to zero 
at all: this was first proved by Koml\'os in 1967~\cite{Ko};
good starting points are Koml\'os' proof 
for $\lim_{d\rightarrow\infty}P_d=O(\frac1{\sqrt{d}})$ given
in \cite[Sect.~XIV.2]{Bo}, and Odlyzko's paper \cite{Odlyzko}.)

Here we only prove the lower bound.
For this, we work in the $\pm1$-model, where
$P_d$ denotes the probability that a random $(d+1)\times(d+1)$-matrix
is singular. In this model, the
probability that two given rows are ``equal or opposite''
is $\frac1{2^d}$, and the same for two given columns.
Altogether there are $2\binom{d+1}2=d^2+d$ 
such events. These are not independent,
but for any two such events the probability that they {\em both\/}
occur is at most $\frac1{2^{2d-1}}$:
if we look at two events that both refer to rows, or both refer
to columns, then the probability that they both occur 
is $\left(\frac1{2^d}\right)^2$;
if we want that two specific rows are equal or opposite, and
two columns are equal or opposite at the same time, then the probability
is $\left(\frac1{2^d}\right)\left(\frac1{2^{d-1}}\right)$.
Thus we may estimate
\[
P_d \ \ \ge\ \ 
(d^2+d) \frac1{2^d}\ \ -\ \ 
\binom{d^2+d}2 \frac1{2^{2d-1}}
\]
and this is larger than $\frac{d^2}{2^d}$ for~$d>10$.
\proofend

It has been conjectured that the lower bound of this
theorem is close to the truth:

\begin{Conjecture}[see~\cite{Odlyzko}, \cite{KKS}]\label{conj:det}
The probability $P_d$ that a random $0/1$-matrix is zero
is dominated by the possibility that one of the rows or
columns is zero, or that two rows are equal, or two columns
are equal:
\[
P_d \ \sim \ 2\binom{d+1}{2}\frac1{2^d}\ \sim\ \frac{d^2}{2^d}.
\]
Equivalently: if a random $\pm1$-matrix of size
$(d+1)\times(d+1)$ is singular, then ``most probably'' two
rows or two columns are equal or opposite. 
\end{Conjecture}


\subsection{High-dimensional random $0/1$-polytopes}

Now we try to describe random $0/1$-polytopes for large~$d$.

\begin{Corollary}\label{cor:komlos}%
With a probability that tends to $1$ for $d\rightarrow\infty$
the following is true:
\begin{itemize}
\item[{\rm(i)}]
Any polynomial number of $0/1$-vectors chosen (independently, with
equal probability) from $\{0,1\}^d$ will be distinct.
\item[{\rm(ii)}]
A set of $d$ randomly chosen $0/1$-points spans a hyperplane
that does not contain the origin~$\zero$.
\item[{\rm(iii)}]
The convex hull of $d+1$ random $0/1$-points
is a $d$-dimensional simplex.
\end{itemize}
\end{Corollary}

\proof
The probability for $n$ random $0/1$-vectors to be distinct is
\[
\left(1-\frac1{2^d}\right)
\left(1-\frac2{2^d}\right)\cdots
\left(1-\frac{n-1}{2^d}\right)
\ \ >\ \ 
\left(1-\frac{n}{2^d}\right)^n
\ \ =\ \ 
\exp\left[n\ln \left(1-\frac{n}{2^d}\right)\right],
\]
and for $n\ll 2^d$ this can be estimated with
$\ln(1-\frac{n}{2^d})\approx-\frac{n}{2^d}$,
so we get a probability of at least $\exp(-\frac{n^2}{2^d})$,
which converges to~$1$ if $\frac{n^2}{2^d}$ tends to zero.

If we choose $d+1$ random points, then by symmetry we may
assume that the first one is~$\zero$. Thus the probability
in question for the third statement, and also for the second one,
is exactly~$1-P_d$, and thus both statements follow from 
Koml\'os' theorem.
\proofend

But one would like to ask more questions. For example:
{\em What is the expected volume of a random simplex?}
It is indeed huge, as one can see from the following observations of 
Szekeres \& Tur\'an \cite{SzekeresTuran}, see Exercise~\ref{Ex:det}:
for a random $0/1$-matrix $A$ of size $d\times d$, the expected value for 
the squared determinant is exactly\footnote{%
The expected value for $\det(A)$ is~$0$ if $d>0$, for symmetry reasons.} 
\[
E(\det(A)^2)\ \ =\ \ \frac{(d+1)!}{2^{2d}}.
\]
But that means that $0/1$-matrices $A$ of determinant 
\[
|\det(A)|\ \ \ge\ \ \frac{\sqrt{(d+1)!}}{2^d}
\]
exist, and are in fact common (``to be expected'').
This is to be compared with the Hadamard upper bound
\[
\det(\Ahat)\ \le\       \sqrt{d+1}^{d+1}, \qquad 
\det( A   )\ \le\ \frac{\sqrt{d+1}^{d+1}}{2^d}.
\]
that we will meet in Section~\ref{sec:0/1-pm1}.

\begin{Proposition}[F\"uredi~\cite{Furedi}]%
For any constant $\varepsilon>0$,
a random $0/1$ polytope with $n\ge(2+\varepsilon)d$ vertices
contains $\frac12\one$, while a random polytope with 
$n\le(2-\varepsilon)d$ vertices does not contain $\frac12\one$,
with probability tending to~$1$ for $d\rightarrow\infty$.
\end{Proposition}

F\"uredi's proof is elementary, combining Koml\'os' theorem
with an estimate about the maximal number of regions
in an arrangement of hyperplanes. Perhaps it can be
adapted to prove that a random $0/1$ polytope 
with $n\ge(2+\varepsilon)d$ vertices should even be centered?

Another question linked to Corollary~\ref{cor:komlos}
is: {\em Can we expect that there will be further
$0/1$-points on the hyperplane spanned by $d$ random points?}
We don't quite know, but the following result points towards an
answer.

\begin{Proposition}[Odlyzko~\cite{Odlyzko}]%
With probability tending to~$1$ for $d\rightarrow\infty$, and
\[
n\ \ \le\ \ d - \frac{10d}{\log d},
\]
$n$ random $0/1$-points span an affine subspace of dimension~$n$ that
does not contain any further $0/1$-point.
\end{Proposition}

One interesting question is whether this result could be extended
to much bigger $n$. Of course, by Corollary
\ref{cor:komlos}(iii) to Koml\'os' theorem the
statement fails (badly) if $n=d+1$, but what about $n=d$?
In other words, is there a high probability that
$d$ random $0/1$-points will span a ``simplex hyperplane''?

Still another, related question is:
{\em If $d$ random points span a hyperplane, is there a
reasonable chance that this hyperplane is very unbalanced, with only
few $0/1$-points on one side?}
This is closely linked
(by ``linearity of expectation'') to the expected number of
facets of a random polytope.

\begin{Proposition}%
There is a constant $c>0$ such that a random
$0/1$-polytope $P\sse[0,1]^d$ with $n\le(1+c)d$ vertices is
``uniform'' in the sense that any $d+1$ points span a $d$-simplex, 
with probability tending to~$1$ for $d\rightarrow\infty$.\\
(In particular, uniform polytopes are simplicial.)
\end{Proposition}

\proof
Let $\const<1$ be a constant such that
$P_d\le\const^d$ holds for all large enough~$d$.
The probability that all $(d+1)$-subsets of a random sequence of $n$
$0/1$-vectors span $d$-simplices is at least 
\[
\mbox{Prob($P$ uniform)}\ \ \ge\ \ 
1\ -\ \binom{n}{d+1} P_d
\ \ >\ \ 
1\ -\ \binom{(1+c)d}{d+1} \const^d
\]
and with $(cd)!\approx\left(\frac{cd}e\right)^{cd}$
we estimate
\[
\binom{(1+c)d}{d+1} \const^d\ \ \le\ \ 
\frac{{\big((1+c)d\big)}}{(cd-1)!}^{cd} \const^d\ \ \approx\ \ 
\left(\frac{e(1+c)}c\right)^{cd} \const^d.
\]
Thus \mbox{\,Prob($P$ uniform)\,} will tend to $1$ 
for large $d$ if
\[
\left(\frac{e(1+c)}c\right)^c\ \ <\ \ \frac1\const.
\]
Thus by Theorem~\ref{thm:komlos} one can take $c=0.00009$.
However, if Conjecture~\ref{conj:det} were true, then one could
indeed take $c=0.27$.
\proofend

Note that if $P$ is simplicial, then $P_1^s$ is a simplex of dimension
at most $d-1$ for $s=0,1$, and thus $P$ has not more than $2d$ vertices.
And simplicial polytopes with $2d$ vertices do indeed exist: but
the only examples that we know are centrally-symmetric cross polytopes, 
which one gets as 
\[
\conv\{\vv_1,\vv_2,\ldots,\vv_d,\one-\vv_1,\one-\vv_2,\ldots,\one-\vv_d\},
\]
where $\vv_1,\vv_2,\ldots,\vv_d\in\{0,1\}^d$ are $d$
affinely independent points whose last coordinate is~$0$.
Are there any other examples? This is not clear, but one may
note that  if $P$ is a $d$-dimensional cross polytope, then it must be
centrally symmetric. In fact, if $\vv,\ww$ are vertices of $P$
that are not adjacent, then they are not both contained in any trivial face
$P_i^s$ (since these faces are simplices), hence they are opposite to
each other in the $d$-cube. 
But is every simplicial $d$-dimensional $0/1$-polytope
with $2d$ vertices necessarily a cross polytope?

%


\section{Cut Polytopes}\label{sec:cut}

The ``special'' $0/1$-polytopes studied in combinatorial optimization
exhibit enormous complexity. One well-studied instance is that of
the symmetric and asymmetric travelling salesman (TSP) polytopes (see
\cite{GroetschelPadberg}), for which 
Billera and Sarangarajan \cite{BS} have
recently shown that \emph{all} $0/1$-polytopes appear as faces.

In this lecture, we discuss basic properties of a different
family of $0/1$-polytopes, the cut polytopes, and of the 
correlation polytopes (a.k.a.\ boolean quadric
polytopes), which are affinely equivalent to them. For all of this
and much more, Deza \& Laurent \cite{DL} provides an excellent and
comprehensive reference.


\subsection{``Small'' cut polytopes}

Let's start with a ``construction by example'' of the ``very small''
cut polytopes; the general prescription will come in the
next section.

A {\em cut\/} in a graph is any edge set of the
form $E(S,V\sm S)=E(V\sm S,S)$, for $S\sse V$.
That is, a cut consists of all edges
that connect a node in~$S$ to a node not in~$S$.
For example, the complete graph~$K_3$ 
\[
\begin{picture}(0,0)%
\epsfig{file=EPS/k3.pstex}%
\end{picture}%
\setlength{\unitlength}{2960sp}%
\begingroup\makeatletter\ifx\SetFigFont\undefined%
\gdef\SetFigFont#1#2#3#4#5{%
  \reset@font\fontsize{#1}{#2pt}%
  \fontfamily{#3}\fontseries{#4}\fontshape{#5}%
  \selectfont}%
\fi\endgroup%
\begin{picture}(1477,1290)(2124,-1936)
\put(2626,-1336){\makebox(0,0)[rb]{\smash{\SetFigFont{12}{14.4}{\rmdefault}{\mddefault}{\updefault}$13$}}}
\put(3376,-1336){\makebox(0,0)[lb]{\smash{\SetFigFont{12}{14.4}{\rmdefault}{\mddefault}{\updefault}$12$}}}
\put(3001,-1936){\makebox(0,0)[b]{\smash{\SetFigFont{12}{14.4}{\rmdefault}{\mddefault}{\updefault}$23$}}}
\put(3001,-811){\makebox(0,0)[b]{\smash{\SetFigFont{9}{10.8}{\rmdefault}{\mddefault}{\updefault}$1$}}}
\put(3601,-1786){\makebox(0,0)[lb]{\smash{\SetFigFont{9}{10.8}{\rmdefault}{\mddefault}{\updefault}$2$}}}
\put(2401,-1786){\makebox(0,0)[rb]{\smash{\SetFigFont{9}{10.8}{\rmdefault}{\mddefault}{\updefault}$3$}}}
\end{picture}

\]
has four cuts:
all the edge sets of size~$2$, as well as the empty set of
edges:
\[
\epsffile{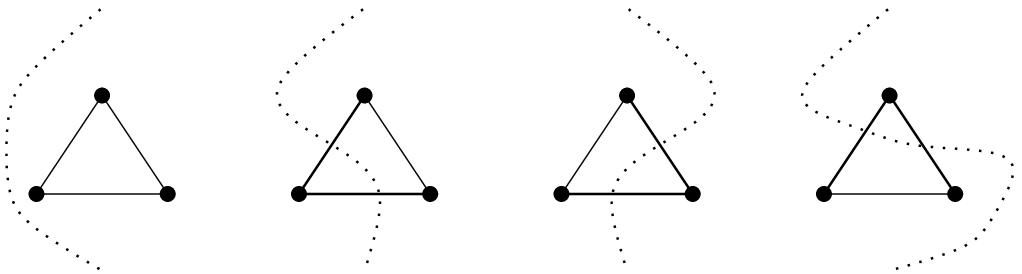}
\]
These cuts can be encoded by their {\em cut vectors\/}
\[
\left(\begin{array}{c}x_{12}\\x_{13}\\x_{23}\end{array}\right)\ \in \{0,1\}^3,
\]
where the $ij$-coordinate records whether
the edge $ij$ is in the cut or not.
The cut polytope is the convex hull of all these cut vectors. 
So, for $K_3$ we get
the cut polytope {\tt CUT3:3-4.poly} as
\[
\CUT(3)\ =\ \conv\Big\{\quad
\left(\begin{array}{c}0\\0\\0\end{array}\right)\quad
\left(\begin{array}{c}0\\1\\1\end{array}\right)\quad
\left(\begin{array}{c}1\\0\\1\end{array}\right)\quad
\left(\begin{array}{c}1\\1\\0\end{array}\right)
\quad\Big\}.
\]
This $0/1$-polytope is the convex hull of all $0/1$-vectors
of even weight (those just happen to be the cuts),
so it is the regular simplex of side-length $\sqrt2$. 
Not a very interesting $0/1$-polytope.

The complete graph $K_4$ has $\binom42=6$ edges,
and altogether $8$ cuts: the empty cut, the
four cuts of size $3$ that separate one vertex from the three others,
and three cuts of size~$4$ that separate two vertices from the two
others. Each cut yields a cut vector 
\[
(x_{12},x_{13},x_{14},x_{23},x_{24},x_{34})^t \ \in \{0,1\}^6.
\]
The resulting polytope {\tt CUT4:6-8.poly}
again has a very simple structure:
it is a sum of two simplices,
\[
\CUT(4)\ \ \cong\ \ \CUT(3)\ *\ \CUT(3)\ \ \cong\ \ \Delta_3 \ *\ \Delta_3.
\]
To see this, note that four of the eight cuts of $K_4$
\[
\epsffile{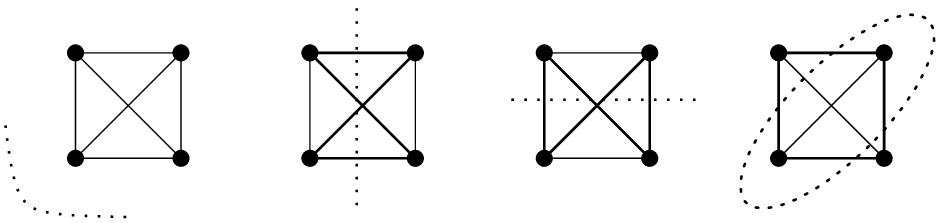}
\]
contain ``none or both'' from each pair of disjoint
edges in~$K_4$, that is,
\[
x_{12}=x_{34}, \qquad
x_{13}=x_{24}, \qquad
x_{14}=x_{23},
\]
so they lie in the $3$-dimensional subspace $U_1$ of $[0,1]^6\sse\R^6$
that is given by these three equations. The other four cuts (of size~$3$)
\[
\epsffile{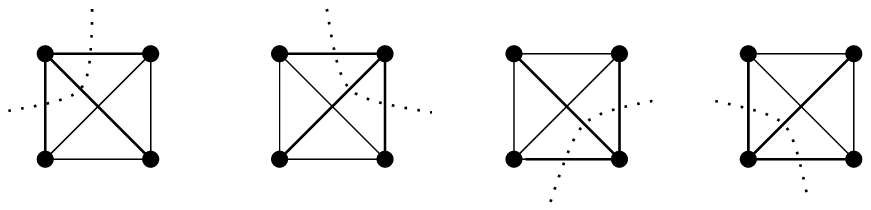}
\]
all contain exactly one edge from each
disjoint pair, that is, they lie in the $3$-dimensional subspace $U_2$
given by 
\[
x_{12}+x_{34}=1, \qquad
x_{13}+x_{24}=1, \qquad
x_{14}+x_{23}=1
\]
and give a $3$-simplex that is equivalent to $\CUT(3)$ in this subspace.
Now $U_1\cap U_2=\{\frac12\one\}$ completes the analysis:
we {\em understand the structure\/} of~$\CUT(4)$.
(Combinatorially, $\CUT(4)$ may also be identified with the cyclic
polytope $C_6(8)$; in particular, it is simplicial, and
neighborly. But nevertheless,
it is not a very interesting polytope.)

And so on \ldots ? It turns out that 
the cut polytopes are much more complicated (``interesting'')
than one might think.


\subsection{Cut polytopes and correlation polytopes}

The definition/construction of the general cut polytopes
follows a general method that has proved to be extremely successful
in combinatorial optimization: The cuts in a complete graph $K_n$
are encoded into the $0/1$-polytope given by their characteristic vectors.

\begin{Definition}[Cut polytopes]%
With every subset $S\subseteq [n]:=\{1,\ldots,n\}$, associate a
$0/1$-vector
\[
\delta(S)\in\{0,1\}^d,\quad d=\binom{n}{2},
\]
by setting (for $1\le i<j\le n$)
\[
\delta(S)_{ij}:=
\left\{
\begin{array}{ll}
1 & \textrm{if }|S\cap\{i,j\}|=1,\\
0 & \textrm{otherwise.}
\end{array}
\right.
\]
Thus we can identify the coordinates $x_{ij}$ of $\R^d$ with the edge set
of $K_n$ (a complete graph with vertex set $[n]$), and the vector
$\delta(S)$ represents the set $\{ij:x_{ij}=1\}$ of edges $ij$ of~$K_n$ 
that connect a vertex in
$S$ with a vertex in $\overline S:=[n]\sm S$, that is, a {\em cut\/}
$E(S,\overline S)$ in $K_n$.

The \emph{cut polytope} $\CUT(n)$ is defined by
\[
\CUT(n)\ \ :=\ \ \conv\,\big\{\,\delta(S):S\subseteq[n]\,\big\}
\quad\subseteq\ \R^d.
\]
\end{Definition}

\begin{Lemma}\label{lemma:cutcentered}%
For every $n\ge 1$, and $d=\binom{n}2$, the cut polytope  $\CUT(n)$ is a
centered $d$-dimensional polytope with $2^{n-1}$ vertices.
\end{Lemma}

\proof
The two sets $S$ and $\overline S$ determine the same cut
$\delta(S)=\delta(\overline S)$, but any two
subsets $S,S'\subseteq[n-1]$ with
$S\not= S'$ determine different cuts $\delta(S)\not=\delta(S')$,
since $S=\{i\in[n-1]:\delta(S)_{in}=1\}$. Thus
\[
\CUT(n)\ \ =\ \ \conv\,\big\{\,\delta(S):S\subseteq[n-1]\,\big\}
\quad\subseteq\ \R^d
\]
has $2^{n-1}$ vertices (corresponding to the $2^{n-1}$ cuts of $K_n$).
If $\CUT(n)\subseteq\R^d$ were not full-dimensional, then it
would satisfy some linear equation:
\[
\aaa^t\xx\ =\ \sum\limits_{i,j}a_{ij}x_{ij}\ =\ \beta\quad\textrm{ for all
}x\in\CUT(n)
\]
for some non-zero $\aaa\in\R^d$. However, the zero cut
$\delta(\emptyset)=\zero\in\CUT(n)$ yields $\beta=0$.
Furthermore, we derive from the sketch below that
\[
\delta(\{i\})+\delta(\{j\})-\delta(\{i,j\})\ \ =\ \ 2\ee_{ij},
\]
\[
\begin{picture}(0,0)%
\epsfig{file=EPS/cut1.pstex}%
\end{picture}%
\setlength{\unitlength}{1973sp}%
\begingroup\makeatletter\ifx\SetFigFont\undefined%
\gdef\SetFigFont#1#2#3#4#5{%
  \reset@font\fontsize{#1}{#2pt}%
  \fontfamily{#3}\fontseries{#4}\fontshape{#5}%
  \selectfont}%
\fi\endgroup%
\begin{picture}(2270,1898)(1714,-2544)
\put(2176,-961){\makebox(0,0)[rb]{\smash{\SetFigFont{12}{14.4}{\rmdefault}{\mddefault}{\updefault}$i$}}}
\put(3601,-961){\makebox(0,0)[lb]{\smash{\SetFigFont{12}{14.4}{\rmdefault}{\mddefault}{\updefault}$j$}}}
\end{picture}

\]
so $\aaa^t\delta(S)=0$ for all $S\subseteq[n]$ implies that 
\[
\aaa^t(2\ee_{ij})\ =\ 2a_{ij}\ =\ 0\quad\textrm{ for all }\{i,j\}\subseteq[n],
\]
and thus $\aaa=\zero$.

To see that the cut polytopes are centered, it suffices to
note any edge $ij$ will be contained in a random cut with probability
exactly $\frac12$. Thus the average over all vertices of
$\CUT(n)$ (that is, the centroid
of the set of vertices) is $\frac12\one$.
\proofend

We note one more feature of the polytope $\CUT(n)$: it is very symmetric,
with a vertex-transitive symmetry group.
In fact, every symmetric difference of two cuts is a cut: this
follows from the equation
\[
E(S,\overline S)\ \symmdiff\ E(T,\overline T)\ \ =\ \ 
E(S\symmdiff T,\overline{S\symmdiff T}),
\]
which is best verified and 
visualized in a little picture such as the following:
\[
\begin{picture}(0,0)%
\epsfig{file=EPS/cut2.pstex}%
\end{picture}%
\setlength{\unitlength}{1973sp}%
\begingroup\makeatletter\ifx\SetFigFont\undefined%
\gdef\SetFigFont#1#2#3#4#5{%
  \reset@font\fontsize{#1}{#2pt}%
  \fontfamily{#3}\fontseries{#4}\fontshape{#5}%
  \selectfont}%
\fi\endgroup%
\begin{picture}(3292,3112)(1681,-4208)
\put(2401,-1411){\makebox(0,0)[lb]{\smash{\SetFigFont{12}{14.4}{\rmdefault}{\mddefault}{\updefault}$T$}}}
\put(2251,-2311){\makebox(0,0)[rb]{\smash{\SetFigFont{12}{14.4}{\rmdefault}{\mddefault}{\updefault}$S$}}}
\end{picture}

\]
Thus for any $S\subseteq[n]$ the {\em switching map}
\[
\begin{array}{ll}
\sigma_S: & \R^d\rightarrow \R^d\\
          & x_{ij}\mapsto
\left\{
\begin{array}{ll}
1-x_{ij} & \textrm{if }ij\in E(S,\overline
S),\textrm{i.~e., if }\delta(S)_{ij}=1,\\
x_{ij}    & \textrm{otherwise},
\end{array}
\right.
\end{array}
\]
defines an automorphism of $\CUT(n)$ that takes 
$\delta(T)$ to $\delta(T\symmdiff S)$, and thus takes
the vertex
$\delta(S)$ to the vertex $\delta(\emptyset)=\zero$, and conversely.
Thus, under such switching operations all vertices of $\CUT(n)$
are equivalent!



Next we will look at a different class of important $0/1$-polytopes:
the cut polytopes in (thin) disguise.

\begin{Definition}[Correlation polytopes]%
The $n$-th \emph{correlation polytope} is the convex hull of all
$n\times n$
$0/1$-matrices of rank~$1$:
\[
\COR(n)\ :=\ \conv\{\xx\xx^t:\xx\in\{0,1\}^n\}\ \subseteq\ \R^{n^2}.
\]
\end{Definition}

It is not so hard to see directly that $\COR(n)$ is a polytope of
dimension $\binom{n+1}{2}$ with $2^{n}$ vertices, but the following
observation yields even more.

\begin{Lemma}[de Simone \cite{deSimone}]\label{lemma:corcut}%
For $n\ge 2$ and $d=\binom{n}2$, there is a linear map
\[
\gamma:\ \R^{(n-1)^2}\ \ \longrightarrow\ \ \R^d
\]
that induces an affine isomorphism of polytopes
\[
\gamma:\ \COR(n-1)\ \ \cong\ \ \CUT(n).
\]
\end{Lemma}

\proof
For every correlation matrix $\xx\xx^t$ we can extract the vector
$\xx$ from its diagonal, from this derive a set
$S\subxx:=\{i\in[n-1]:x_i=1\}$, and thus get the cut vector
$\delta(S\subxx)$. Furthermore, the components of $\delta(S\subxx)$
can be derived as {\em linear\/} combinations of the entries 
of~$\xx\xx^t$:
\begin{eqnarray*}
\delta_{ij}&:=&x_i(1-x_j)+x_j(1-x_i)=x_{ii}-x_{ij}+x_{jj}-x_{ji},\\
\textrm{ and }\quad
\delta_{in}&:=&x_i=x_{ii}.
\end{eqnarray*}
This defines a linear map
$\gamma:\R^{(n-1)^2}\rightarrow\R^d$ which maps correlation matrices
to cut vectors: $\gamma(\xx\xx^t)=\delta(S\subxx)$, and thus
$\gamma(\COR(n-1))=\CUT(n)$.
An inverse map is obtained by taking
\begin{eqnarray*}
x_{ii}&:=&\delta_{in}\\
\textrm{ and }\quad
x_{ij}=x_{ji}&:=&\textstyle{\frac{1}{2}}(x_{ii}+x_{jj}-\delta_{ij})=
\textstyle{\frac{1}{2}}(\delta_{in}+\delta_{jn}-\delta_{ij}).
\end{eqnarray*}
The image of this inverse map consists of only symmetric matrices in
$\R^{(n-1)^2}$,
which describes the $\binom{n}{2}$-dimensional subspace of
$\R^{(n-1)^2}$
that is spanned by~$\COR(n-1)$.
\proofend

Note that the isomorphism of Lemma~\ref{lemma:corcut} is
{\em not\/} a $0/1$-equivalence --- in fact the polytopes are
not $0/1$-equivalent, even in their full-dimensional versions. 
For example cut polytopes are centered (Lemma~\ref{lemma:cutcentered}),
but the correlation polytopes are not:
$\COR(n)$ contains the point $\frac12\one=\frac12(\zero+\one)$, but 
this point lies in the boundary, since
$x_{11}\ge x_{12}$ is valid for all vertices of
$\COR(n)$, and not for all of them with equality.
 
We now record a remarkable property of the correlation polytopes
(and of cut polytopes, via Lemma~\ref{lemma:corcut}): 

\begin{Proposition}\label{Pr:3-neighborly}%
Any three vertices of $\COR(n)$ determine a triangular face
$F\cong\Delta_2$, that is,
$\COR(n)$ is $3$-neighborly.
\end{Proposition}

\proof
Using the symmetry of $\CUT(n+1)$, and its affine equivalence with
$\COR(n)$, we may assume that one of the three vertices of $\COR(n)$ is
$\zero\zero^t$, while the others are $\uu\uu^t$ and $\vv\vv^t$. The
vectors
$\uu,\vv\in\R^n$ span a $2$-dimensional subspace $U\subseteq \R^n$,
which may or may not contain a fourth $0/1$-vector $\yy\in\R^n$,
but no fifth vector.
However, if
there is such a fourth vector $\yy$, then we may assume that
$\yy=\uu+\vv$
(possibly after exchanging $\yy$ with $\uu$ or with $\vv$).

Now take a \emph{generic} vector $\hh\in\R^n$ that is orthogonal
to $U$ --- such a vector will satisfy $\hh^t\uu=\hh^t\vv=\hh^t\zero=0$, and
also $\hh^t\yy=0$ if $\yy$ exists, but $\hh^t\xx\not= 0$ for any other
$\xx\in\{0,1\}^n$.
Then a little computation shows that the standard scalar product on
$\R^{n^2}$ with $\hh\hh^t$ defines a \emph{linear} function on
$\COR(n)$ that is minimized by $\zero\zero^t,\uu\uu^t,\vv\vv^t$,
and by $\yy\yy^t$ if this $\yy$ exists, but by no other vertex of
$\COR(n)$:
\[
\big\langle\hh\hh^t,\xx\xx^t\big\rangle
\ \ =\ \ \sum\limits_{1\le i,j\le n}(\hh\hh^t)_{ij}(\xx\xx^t)_{ij}
\ \ =\ \ \sum\limits_{1\le i\le n}
         \sum\limits_{1\le j\le n} h_i h_j x_i x_j\ \ =
\]
\[
=\ \ \Big(\sum\limits_{1\le i\le n}h_i x_i\Big)
     \Big(\sum\limits_{1\le j\le n}h_j x_j\Big)
\ \ =\ \ (\hh^t\xx)^2\ \ \ge\ 0.
\]
Now if there is no ``fourth man'' $\yy$, then this proves that
$\conv(\{\uu\uu^t,\vv\vv^t,\zero\zero^t\})$ is a (triangular) face of
$\COR(n)$. If, however, $\yy=\uu+\vv$ is present (that is,
$\uu+\vv\in\{0,1\}^n$, and thus $\uu^t\vv=0$), then we obtain that 
\[
F\ \ :=\ \ \conv\,\big(\,
\{\zero\zero^t,\uu\uu^t,\vv\vv^t,(\uu+\vv)(\uu+\vv)^t\}\,\big)
\]
is a face of $\COR(n)$. 
We have to show that this face is a tetrahedron, not a $2$-face.

Since 
$\uu^t\vv=0$ with $\uu,\vv\neq\zero$, we can take indices $i,j$
with $u_i=1$, $v_i=0$ and $v_j=1$, $u_j=0$, so that 
\[
(\uu\uu^t)_{ij}=u_iu_j=0,\qquad
(\vv\vv^t)_{ij}=v_iv_j=0,\qquad
(\yy\yy^t)_{ij}=(u_i+v_i)(u_j+v_j)=1.
\]
Thus $\yy\yy^t$ cannot be linearly dependent of $\uu\uu^t$ and $\vv\vv^t$,
and it is also clear that $\uu\uu^t$ and $\vv\vv^t$
are distinct $0/1$-vectors and hence linearly independent.
Thus $\uu\uu^t,\vv\vv^t$ and $(\uu+\vv)(\uu+\vv)^t$ are
linearly independent, and hence $F$ is a tetrahedron face of~$\COR(n)$.
\proofend

This result is best possible, since 
$\CUT(n)$ is not $4$-neighborly in general: for
this we note (for $n\ge 3$) that
\[
\delta(\emptyset)+\delta(\{1,2\})+\delta(\{1,3\})+\delta(\{2,3\})=
\delta(\{1\})+\delta(\{2\})+\delta(\{3\})+\delta(\{1,2,3\})
\]
which implies that the four vectors on either side of the equation 
(which are distinct vectors for $n\ge4$) do
\emph{not} form a tetrahedron face of $\CUT(n)$.

Proposition~\ref{Pr:3-neighborly} implies that
$\CUT(n)$ is {\em $5$-simplicial}, that is, all
the $5$-dimensional faces of $\CUT(n)$ are simplices
(Exercise \ref{Ex:6-simplicial}).
On the other hand, the cut polytopes are not
$6$-simplicial: $\CUT(4)$ is $6$-dimensional,
but it is not a simplex.
(Check {\tt SIMPLICIALITY} for the cut polytopes in the
{\tt polymake} database!


\begin{Corollary}\label{cor:neigh}%
For every dimension $d=\binom{n}{2}$, there is a $3$-neighborly
$0/1$-polytope with more than $2^{\sqrt{2d}-1/2}$ vertices.
\end{Corollary}

\proof
Take $\CUT(n)$, whose number of vertices is $2^{n-1}$, with
$n=\frac{1}{2}+\sqrt{2d+\frac14}$.
\proofend


\subsection{Many facets?}

Here I would also like to give --- at least --- a rough estimate of
the number of facets of $\CUT(n)$ for large $n$, but that seems
not
that easy to get. We note that 
\[
\CUT(2)\ \cong\ \Delta_1\quad\textrm{ and }\quad\CUT(3)\ \cong\ \Delta_3
\]
are simplices, while computation in ``small'' dimensions (see the
{\tt polymake} database and SMAPO~\cite{Ch}) yields

\begin{tabular}{lll}
\CUT (4) & has dimension $d=\hphantom{0}6$ & and $16=1.5874^d$ facets,\\
\CUT (5) & has dimension $d=10$ & and $56=1.4956^d$ facets,\\
\CUT (6) & has dimension $d=15$ & and $368=1.4827^d$ facets,\\
\CUT (7) & has dimension $d=21$ & and $116764=1.7430^d$ facets,\\
\CUT (8) & has dimension $d=28$ & and $217093472=1.9849^d$ facets,\\
\CUT (9) & has dimension $d=36$ & and at least $12246651158320=2.3097^d$
facets.
\end{tabular}

This suggests that $\CUT(n)$ has more than $d^{cd}$ facets,
for some $c>0$: prove this!


\section{The Size of Coefficients}\label{sec:coeff}

Gr\"otschel, Lov\'asz \& Schrijver \cite{GLS}, 
in their study of the
ellipsoid method and its (fundamental) role in optimization,
introduced
the notion of the \emph{facet complexity\/} of a polyhedron. This is
roughly
the maximal number of bits that is necessary to represent one single
facet by an inequality (with rational coefficients). They showed that
for polyhedra with bounded facet complexity, optimization and
separation are equivalent. Thus, the complexity of the facets is more
important in this context than the number of facets. The following will
imply that the facet complexity of an $n$-dimensional $0/1$-polytope is
$O(n^2\log n)$: this is a polynomial bound, and thus ``good enough'' for the
ellipsoid method.

The question about the maximal facet complexity of $0/1$-polytopes
can also be phrased differently: it asks {\em How large integers
(rationals) may occur in the $\cal H$-presentation of a $0/1$-polytope?\/}
The bad news is that the integer
coefficients
that appear in the inequality description of a $0/1$-polytope may be
\emph{huge}. This is ``bad'': it means that 
all kinds of algorithms, from cutting plane procedures to 
convex hull algorithms --- used to compute the facets
of a given polytope --- are threatened by ``integer
overflow'' even in the case of $0/1$-polytopes.

The main source for this lecture is a recent paper by Noga Alon
and V{\v a}n H.~V{\~u} \cite{AV}, which rests on a construction of
Johan H{\aa}stad \cite{Hastad} from 1992.


\subsection{Experimental evidence}

What do we mean by the size of the coefficients of the
facets? For this we write each facet-defining inequality
of a full-dimensional (!) $0/1$-polytope uniquely in the
normal form 
\[
\pm c_0\, \pm c_1x_1 \pm c_2x_2 \pm \ldots \pm c_dx_d\ \ \ge\ 0,
\]
for non-negative integers $c_0,c_1,\ldots,c_d$
with greatest common divisor~$1$.
By the {\em greatest coefficient\/} we mean
$\max\{c_1,\ldots,c_d\}$.
For example, the inequality
\[
19 -12x_1 -18x_2 -3x_3 -1x_4 + 10x_4 -11x_6 + 4x_7 -5x_8\ \ \ge\ \ 0
\]
from {\tt CF:8-9.poly} has greatest coefficient $18$.

The concept of greatest coefficient is invariant under permuting coordinates
(obviously), but also under switching (the substitution
$x_i\leftrightarrow 1-x_i$ just switches the sign in front of $c_ix_i$,
but not the size of the coefficient). It also changes the constant
coefficient $c_0$, but we ignore these anyway. Note that for 
$0/1$-polytopes we always have
$c_0\le c_1+\ldots + c_d$, since the facet-defining 
inequality must be satisfied by some $0/1$-point with equality. 
We will, however, apply the concept of ``greatest coefficients''
only in the case of full-dimensional polytopes,
since otherwise the ``defining inequality of a facet'' is not unique,
which makes things more complicated.
\newcommand{\coeff}{\operatorname{coeff}}

With these precautions, we can look up the largest coefficient
$\coeff(d)$
that appears in a facet-defining inequality for a $d$-dimensional
$0/1$-polytope, and for low dimensions $d$ we find the following:
\[
\begin{array}{r|l@{\quad\quad}l}
d & \ \ \ \coeff(d) & \mbox{example}\\ 
\hline
\\[-4mm]
3 & \qquad= \hphantom{0}1\\
4 & \qquad= \hphantom{0}2 & \mbox{\tt CF:4-5.poly} \\
5 & \qquad= \hphantom{0}3 & \mbox{\tt CF:5-6.poly} \\
6 & \qquad= \hphantom{0}5 & \mbox{\tt CF:6-7.poly} \\
7 & \qquad= \hphantom{0}9 & \mbox{\tt CF:7-8.poly} \\
8 & \qquad= 18  & \mbox{\tt CF:8-9.poly}\\
9 & \qquad\ge 42  & \mbox{\tt CF:9-10.poly}\\
10 & \qquad\ge 96  & \mbox{\tt CF:10-11.poly}
\end{array}
\]
Here the values for $d\le 8$ are from complete enumeration,
the values for $d>8$ were taken from Aichholzer~\cite[p.~111]{Ai-dipl}.
The data for $d\le10$ do not, however, provide enough evidence
to guess the truth.


\subsection{The Alon-V\~u theorem and some applications}\label{sec:0/1-pm1}

Let $A$ be a $0/1$-matrix of size $n\times n$.
The question {\em How bad can $A$ be?\/} has many aspects.
Here we will first look (again) at the maximal size of a determinant $\det(A)$.
Then we get to the Alon-V\~u theorem about 
the maximal size of entries of~$A^{-1}$, and to its consequences
for the arithmetics (large coefficients) and the geometry 
(e.~g.\ flatness) of $0/1$-polytopes.

Denote by $\rho_n$ the maximal determinant of a 
$0/1$-matrix of size $n\times n$.
The exact value of $\rho_n$ seems to be known for all $n<18$,
except for $n=14$, where the following table quotes a conjecture.

\[
\begin{array}{r|rl}
n & \rho_n \\
\hline
1  & 1 \\
2  & 1 \\
3  & 2 \\
4  & 3 \\
5  & 5 \\
6  & 9 \\
7  & 32 \\
8  & 56 \\
9  & 144 \\
10 & 320 \\
11 & 1458 \\
12 & 3645 \\
13 & 9477 \\
14 & 25515  & \mbox{\ \ (?, Smith~\cite{Smith87}, Cohn~\cite{Cohn})}\\
15 & 131072 \\
16 & 327680 &  \ldots
\end{array}
\]
Matrices that achieve these values may be obtained from a
web page by Dowdeswell, Neubauer, Solomon \& Tumer \cite{DowdeswellSolomon}.

\begin{Lemma}[The Hadamard bound]
The maximal determinant of a $0/1$-matrix of size $n\times n$ is bounded by
\[
\rho_n\ \ \le\ \ 2\left(\frac{\sqrt{n+1}}2\right)^{n+1}.
\]
\end{Lemma}

\proof
The Hadamard inequality states that the determinant 
of a square matrix is at most the product of the lengths
of its columns, with equality (in the nonsingular case) if and only if
all columns are orthogonal to each other.
Applied to the case of a $\pm1$-matrix $\Ahat$ of size $(n+1)\times(n+1)$,
this yields  
\begin{equation}\tag{$*$}
\det(\Ahat)\ \ \le\ \ \sqrt{n+1}^{~n+1}.
\end{equation}
We transfer this result to $n\times n$ $0/1$-matrices $A$ via
Proposition~\ref{Prop:biject}, and get
\[
\det(A)\ \ \le\ \ \frac{\sqrt{n+1}^{~n+1}}{2^n},
\]
as claimed.
\proofend

A matrix $\Ahat\in\{-1,+1\}^{(n+1)\times(n+1)}$
that achieves equality in ($*$) is known
as a {\em Hadamard matrix}. It is not hard to show that
for this a condition is that $n+1$ is $1,2$ or a multiple of~$4$.
It is conjectured that these conditions are also 
sufficient, but for many values $n+1\ge428$ this is not known.
We refer to Hudelson, Klee \& Larman~\cite{HKL} 
for an extensive, recent survey with pointers to the vast literature 
related to the Hadamard determinant problem.
For the cases where $n+1$ is not a multiple of $4$ one has
slightly better estimates (by a constant factor) than the Hadamard bound;
see Neubauer and Radcliffe \cite{NeubauerRadcliffe}.
Certainly for our purposes we may consider the Hadamard bound as
``essentially sharp.''

Now assume additionally that $A\in\{0,1\}^{n\times n}$ 
is invertible (of determinant
$\det(A)\neq0$), consider the inverse $B:=A^{-1}$, and let
\[
\chi(A):=\max\limits_{1\le i,j\le n}|b_{ij}|=\|B\|_{\infty}
\]
the largest absolute value of an entry of $A^{-1}$. These entries are
---
by Cramer's rule --- given by
\[
b_{ij}=(-1)^{i+j}\det (A_{ij})/\det(A),
\]
where $A_{ij}$ is obtained from $A$ by deleting the $i$-th row and
the
$j$-th column. Let $\chi(n)$ denote the maximal entry in the
inverse of any invertible $0/1$-matrix of size $n\times n$.

\begin{Theorem}[Alon \& V{\~u} \cite{AV}]%
The maximal absolute value of an entry in the inverse of an invertible
$0/1$-matrix of size $n\times n$ can be bounded by
\[
\frac{n^{n/2}}{2^{2n+o(n)}}
\ \ \le\ \ \chi(n)
\ \ \le\ \ \rho_{n-1}
\ \ \le\ \ 
\frac{n^{n/2}}{2^{n-1}}.
\]
Furthermore, $0/1$-matrices that realize the lower bound can be
effectively constructed. (An even better lower bound, by a factor
of $2^n$, is achieved in the case where $n$ is a power of~$2$.)
\end{Theorem}



Before we look at the proof of this theorem, we derive two (quite
immediate) applications to the geometry of $0/1$-polytopes.
First, let as above  $\coeff(d)$ denote the largest $c_i$
that can appear in a reduced inequality 
\[
\pm c_0\, \pm c_1x_1 \pm c_2x_2 \pm \ldots \pm c_dx_d\ \ \ge\ 0,
\]
that defines a facet of a $d$-dimensional $0/1$-polytope in $\R^d$.
(Here the $c_i$ are non-negative integers, with 
$\gcd(c_1,\ldots,c_d)=1$; by switching, we may assume that
$c_0=0$ if we want~to.)

\begin{Corollary}[Huge coefficients~\cite{AV}]%
The largest integer coefficient $\coeff(d)$ in
the facet description of a full-dimensional $0/1$-polytope in~$\R^d$
satisfies
\[
\frac{(d-1)^{(d-1)/2}}{2^{2d+o(d)}}
\ \ \le\ \ \chi(d-1)
\ \ \le\ \ \coeff(d)
\ \ \le\ \ \rho_{d-1}
\ \ \le\ \  \frac{d^{d/2}}{2^{d-1}}.
\]
\end{Corollary}

\proof
Let $\{\zero,\vv_1,\ldots,\vv_{d-1}\}\sse\{0,1\}^d$ be points that span
a hyperplane $H$ in $\R^d$, and let
$V=(\vv_1,\ldots,\vv_{d-1})^t\in\{0,1\}^{(d-1)\times d}$. Then an
equation that defines $H$ is given by 
$\cc^t x=0$, with
$c_i=\pm\det(V_i)$, where $V_i\in\{0,1\}^{(d-1)\times(d-1)}$ is obtained
from $V$ by deleting the $i$th column.
Thus we get the upper bound $\coeff(d)\le\rho_{d-1}$ by definition.

For the lower bound $\chi(d-1)\le\coeff(d)$ we
start with a matrix $A\in\{0,1\}^{(d-1)\times(d-1)}$
such that 
$\chi(A)=|\det A_{11}/\det A|=\chi(d-1)$, and let 
$V:=(A,\ee_1)\in\{0,1\}^{(d-1)\times d}$. 
Then $|\det V_d|=|\det A|$, while $|\det V_1|=|\det A_{11}|$. 
Thus for the coefficients $c_i=\pm\det(V_i)$
of a corresponding inequality $\cc^t \xx\ge 0$ we get
\[
|c_1/c_d|\ =\ |\det A_{11}/\det A_1|\ =\ \chi(d-1),
\]
and thus for any integral inequality which defines a
facet that lies in our hyperplane $H=\{\xx\in\R^d: \cc^t \xx= 0\}$
we have $c_1\ge\chi(d-1)$. 

A simplex for which this $H$ defines a facet is, for example, given by
the convex hull of $\zero$ and $\ee_1$ together with the rows of~$V$.
This simplex has
determinant $\det(A_{11})$, which will be huge for the matrices $A$
constructed for the Alon-V\~u theorem.
\proofend

\newcommand{\flt}{\operatorname{flat}}
\newcommand{\distance}{\operatorname{dist}}
\begin{Corollary}[Flat \boldmath{$0/1$}-simplices~\cite{AV}]%
The minimal positive distance $\flt(d)$ of a $0/1$-vector from a
hyperplane that is spanned by $0/1$-vectors in $\R^d$ satisfies
\[
\frac{2^{d-1}}{\sqrt{d}^{~d+1}}
\ \ \le\ \  \frac{1}{\sqrt{d}\,\rho_{d-1}}
\ \ \le\ \  \flt(d)
\ \ \le\ \  \frac{1}{\chi(d)}
\ \ \le\ \  \left(\frac{1}{d}\right)^{d/2}2^{d(2+o(1))}.
\]
\end{Corollary}

\proof
Let $H=\aff\{\zero , \vv_2,\ldots, \vv_d\}$ be a hyperplane
under consideration (we may assume that it contains the origin) and
let
$\vv_1\in\{0,1\}^d\setminus H$. Then there is
an integral normal vector $\cc$ to~$H$
with $c_i=\pm \det(A_{i1})$, for the square matrix
$A:=(\vv_1,\vv_2,\ldots,\vv_d)\in\{0,1\}^{d\times d}$.
From $\vv_1\not\in H$ we get $|\vv_1^t \cc|\ge 1$, while the length of $\cc$ 
is bounded by 
\[
\|\cc\|\ \ \le\ \ \sqrt{d}\,\|\cc\|_{\infty}\ \ \le\ \ \sqrt{d}\,\rho_{d-1},
\]
and thus
\[
\distance(\vv_1,H)\ =\ 
\frac{|\vv_1^t\cc|}{\|\cc\|}
\ \ \ge\ \ \frac1{\|\cc\|}
\ \ \ge\ \ \frac1{\sqrt{d}\,\rho_{d-1}}.
\]
For the upper bound, take an $A$ that achieves 
$\chi(A)=|\det(A_{11})/\det(A)|=\chi(d)$. Then
\[
\frac1{\chi(d)}
\ =\ \frac{|\det(A)|}{|\det(A_{11})|}
\ =\ \frac{\Vol(\conv\{\zero,\vv_1,\vv_2,\ldots,\vv_d\})}
          {\Vol(\conv\{\zero,\ee_1,\vv_2,\ldots,\vv_d\})}
\ =\   \frac{\distance(\vv_1,H)}
            {\distance(\ee_1,H)}
\ \ge\ \frac{\distance(\vv_1,H)}{1},
\]
where the last ``='' is since we are considering two simplices with a
common facet, and the inequality is from  
$\distance(\ee_1,H)\le \distance(\ee_1,\zero)=1$.
\proofend

\proof 
We now survey the main parts of the proof of the Alon-V\~u theorem,
following \cite{AV}.

{\bf (1) The upper bound. }
For the upper bound $\chi(n)\le\rho_{n-1}$ 
we use that the entries of $A^{-1}$ can be written as
\[
b_{ij}\ \ =\ \ (-1)^{i+j}\frac{\det(A_{ij})}{\det(A)},
\]
where the cofactors $A_{ij}\in\{0,1\}^{(n-1)\times(n-1)}$
satisfy $|\det(A_{ij})|\le\rho_{n-1}$ by definition,
and the invertible matrix $A$ satisfies $|\det(A)|\ge1$
since it is integral. 
%
%


{\bf (2) Super-multiplicativity. }
%
%
For the lower bound it is sufficient to construct ``bad'' matrices 
of size $2^m\times 2^m$, because of the following simple construction,
which establishes
\[
\chi(n_1+n_2)\ \ \ge\ \ \chi(n_1)\cdot\chi(n_2).
\]
Take ``bad'' invertible $0/1$-matrices $A$ and $B$ of sizes
$n_1\times n_1$ and $n_2\times n_2$, such that
$\chi(A)=|\det A_{n_1,n_1}/\det A|$
and $\chi(B)=|\det B_{11}/\det B|$.
Then the matrix
\[
A\diamond B\ \ :=\ \ 
\left(
\begin{array}{cccccc}
~ & ~ & ~ & ~ & ~ & ~\\
~ & A & ~ & ~ & 0 & ~\\
~ & ~ & ~ & ~ & ~ & ~\\
0 & \cdots & 1 & ~ & ~ & ~\\
\vdots & & \vdots & ~ & B & ~\\
0 & \cdots & 0 & ~ & ~ & ~
\end{array}
\right)
\]
has determinant $\det(A\diamond B)=\det(A)\cdot\det(B)$
and the submatrix
\[
(A\diamond B)_{n_1,n_1+1}\ \ =\ \ 
\left(
\begin{array}{ccc@{\quad}cc}
~      & ~           & *      & ~      & ~\\
~      & A_{n_1,n_1} & *      & 0      & ~\\[2mm]
0      & \cdots      & 1      & *      &  *\\
\vdots &             & \vdots & B_{11} & ~\\
0      & \cdots & 0           & ~      & ~
\end{array}
\right)
\]
has determinant $\det A_{n_1,n_1}\det B_{11}$, which establishes
\[
\chi(A\diamond B)\ \ge\ \chi(A)\,\chi(B).
\]
Thus --- modulo an annoying computation that you may find in
\cite[Sect.~2.4]{AV} --- it suffices to establish the lower bound of 
the Alon-V\~u theorem for $n=2^m$.


{\bf (3) The construction. }
Here comes the key part of the proof: an ingenious construction 
of a ``bad'' $\pm1$-matrix whose size is a power of~$2$.
Thus we prove that for $n=2^m$ one can construct an invertible matrix
$A\in\{+1,-1\}^{n\times n}$ with
\[
\chi(A)\ \ =\ \ n^{n/2}\left(\frac{1}{2}\right)^{n+o(n)}
\]
and then use Proposition~\ref{Prop:biject}. 
For this, the following is an explicit recipe. Perhaps you want
to ``do it'' for $m=3$, $n=8$?
\begin{itemize}
\item[(i)] 
Choose an ordering $\alpha_1,\alpha_2,\ldots,\alpha_n$ on the collection
of all $2^m=n$  subsets of $[m]=\{1,2,\ldots,m\}$, such that 
$|\alpha_i|\le|\alpha_{i+1}|$ and $|\alpha_i\symmdiff \alpha_{i+1}|\le 2$
holds for all~$i$. This is not hard to do.
\item[(ii)] 
The matrix $Q\in\{+1,-1\}^{n\times n}$ given by
$q_{ij}:=(-1)^{|\alpha_i\cap\alpha_j|}$ is a symmetric Hadamard
matrix (in
fact, in lexicographic ordering of the rows and columns this is the
``obvious'' Hadamard matrix of order $2^m$). Thus $Q^2=n I_n$,
$Q^{-1}=\frac{1}{n}Q$, and $\det(Q)=n^{n/2}$.
\item[(iii)] 
We construct a lower triangular matrix $L\in\Q^{n\times n}$
row-by-row, with $(1,0,\ldots,0)$ as the first row. For $i>1$ define
$A_i:=\alpha_{i-1}\cup\alpha_i$ and
\[
F_i:=
\left\{
\begin{array}{ll}
\{\alpha_s:\alpha_s\subseteq A_i,\ 
|\alpha_s\cap(\alpha_{i-1}\symmdiff\alpha_i)|=1 & 
                 \textrm{if }|\alpha_{i-1}\symmdiff\alpha_i|=2,\\
\{\alpha_s:\alpha_s\subseteq
A_i=\alpha_i\} & \textrm{if }|\alpha_{i-1}\symmdiff\alpha_i|=1,
\end{array}
\right.
\]
so that both $\alpha_{i-1},\alpha_i\in F_i$ and 
$|F_i|=2^k$ hold in both cases, for
\[
k\ :=\ |\alpha_i|.
\]
Then for $1<i\le n$ and $1\le j\le n$ we set
\[
\ell_{ij}\ \ :=\ \ 
\left\{
\begin{array}{ll}
0 & \textrm{if}\quad \alpha_j\not\in F_i,\\
\left(\frac{1}{2}\right)^{k-1}-1 & \textrm{if }j=i-1,\textrm{ and}
\\
\left(\frac{1}{2}\right)^{k-1} & \textrm{otherwise}.
\end{array}
\right.
\]

\item[(iv)] 
We define $A:= LQ$. A simple computation shows
that $a_{ij}\in\{+1,-1\}$ holds for all~$i,j$.
The determinant of~$A$ is $2^{n-1}$, since
$\det(Q)=n^{n/2}=2^{m2^{m-1}}$ and
\[
\det(L)\ =\ 
\prod\limits^n_{i=1}\ell_{ii}\ =\ 
\prod\limits^m_{k=1}\left(\frac{1}{2}\right)^{(k-1)\binom{m}{k}}\ =\ 
\left(\frac{1}{2}\right)^{\sum\limits^m_{k=1}(k-1)\binom{m}{k}},
\]
with
\[
\sum\limits^m_{k=1}(k-1)\binom{m}{k}\ =\ 
\sum\limits^m_{k=1}k\binom{m}{k}-\sum\limits^m_{k=1}\binom{m}{k}\ =\ 
m2^{m-1}-2^m+1.
\]
Thus $|\det(A)|$ has the minimal possible value for an invertible
$0/1$-matrix of size~$n\times n$.
 
\item[(v)] 
Take $i_0:=2+m+\binom{m}{2}$, which is the smallest index with
$|\alpha_{i_0}|\ge 3$. We solve the system $L\xx=\ee_{i_0}$. This
is easy since $L$ is lower triangular:\\
$x_i=0$ for $i<i_0$,\\
$x_{i_0}=1/\ell_{i_0 i_0}=4$ since $\ell_{i_0
  i_0}=\frac1{2^{3-1}}=\frac14$,\\ 
and for $i>i_0$ we can solve recursively:
\begin{equation}\tag{$*$}
x_i\ \ =\ \ (2^{k-1}-1)x_{i-1}-
\sum\limits_{\alpha_j\in F_i\sm\{\alpha_i,\ \alpha_{i-1}\}}x_j
\qquad\textrm{for }k=|\alpha_i|.
\end{equation}
Using induction, we now verify that the $x_i$ are positive and
\begin{equation}\tag{$**$}
x_i\ \ >\ \ (2^{k-1}-2)x_{i-1}\textrm{\quad for }i>i_0.
\end{equation}
Indeed, this holds for $i=i_0+1$, and by induction
(with $k\ge 3$, so $2^{k-1}-2\ge 2$) we have
\[
x_{i-1}\ >\ 2 x_{i-2}\ >\ 4x_{i-3}\ >\ \ldots
\]
Thus the sum in $(*)$ is smaller than
\[
\frac{1}{2}x_{i-1}+\frac{1}{4}x_{i-1}+\ldots\ \ =\ \ \sum\limits_{t\ge 1}
\frac{1}{2^t}x_{i-1}\ \ =\ \ x_{i-1}.
\]
Using this estimate in $(*)$ we get for $i>i_0$ that 
\begin{equation}\tag{$*{*}*$}
x_i\ \ >\ \ (2^{k-1}-1)x_{i-1}-x_{i-1}\ \ =\ \ (2^{k-1}-2)x_{i-1}.
\end{equation}
Iteration of the recursion $(**)$, with a start at $x_{i_0}>2$, now yields
\[
x_n\ \ >\ \ 
\prod^m_{k=3}(2^{k-1}-2)^{\binom{m}{k}}\ \ =\ \ 
\prod^m_{k=3} 2^{(k-1)\binom{m}{k}}\ 
\prod^m_{k=3}\Big(1-\frac{2}{2^{k-1}}\Big)^{\binom{m}{k}}
\]
where the first product is $2^N$ with
\[
N\ \ =\ \ \sum^m_{k=1}(k-1)\binom{m}{k}-\binom{m}{2}\ \ =\ \ 
                           m2^{m-1}-2^m-\binom{m}{2}
\]
using the same sum as in (iv), and thus 
\[
2^N\ =\ 
2^{m2^{m-1}-2^m-\binom{m}{2}}\ \ =\ \ 
\frac{n^{n/2}}{2^{n+\binom{m}{2}}}\ \ =\ \ 
n^{n/2}\left(\frac{1}{2}\right)^{n+o(n)}.
\]
Now we use that $1-x\ge\frac1{2^{2x}}$ for $0\le x\le \frac12$ and
thus estimate that the second product is at least $(\frac12)^M$ for
\[
M\ =\ 
2\sum_{k=3}^m \frac1{2^{k-2}}\binom{m}{k}\ \ <\ \ 
8\sum_{k=0}^m \frac1{2^k    }\binom{m}{k}\ \ =\ \ 
8\left(\frac32\right)^m\ =\ 
8\,n^{\log3/2}\ =\ o(n).
\]
Taken together, we have verified that
\[
x_n\ \ = \ \ n^{n/2}\left(\frac{1}{2}\right)^{n+o(n)}.
\]
\item[(vi)] 
The rest is easy: to get the $i_0$-th column of $A^{-1}$,
we solve the system
\[
A \yy\ =\ \ee_{i_0}\ \ \Longleftrightarrow\ \  
LQ\yy\ =\ \ee_{i_0}\ \ \Longleftrightarrow\ \  
 Q\yy\ =\ \xx\textrm{ and }L\xx\ =\ \ee_{i_0}.
\]
But $Q\yy=\xx$ is easy to solve because of $Q^{-1}=\frac{1}{n}Q$. Thus we
obtain
\[
B_{ii_0}\ \ =\ \ y_i\ \ =\ \ \frac{1}{n}\sum^n_{j=1}q_{ij}x_j.
\]
Here $|q_{ij}|=1$ by construction and
from $(*{*}*)$, for $k\ge4$ ($n\ge16$), we have
\[
x_n\ >\ 4x_{n-1}\ >\ 8x_{n-2}\ >\ \ldots
\]
which yields
\[
B_{ii_0}\ \ =\ \ y_i\ \ >\ \ \frac{1}{n}\left(\frac{1}{2}x_n\right)\ \ =\ \ 
            \frac{1}{2n}x_n\ \ \ge\ \ 
n^{n/2}\left(\frac{1}{2}\right)^{n+o(n)}.
\]  
Thus \emph{all} entries of the $i_0$-column of $A^{-1}$ are ``huge.''
\proofend
\end{itemize}


\subsection{More experimental evidence}

The Alon-V\~u construction is completely explicit;
you will find corresponding simplices (generated by Michael Joswig)
as {\tt MJ:16-17.poly} and as {\tt MJ:32-33.poly}
in the {\tt polymake} database.
The first one is a $16$-dimensional simplex 
with ``$-451$'' appearing as a coefficient.
The second one has dimension~$32$, and here
you'll find tons of coefficients like
``$4964768222$'' that are indeed large enough to
cause trouble for any conventional single-precision 
arithmetic system~\ldots


\section{Further Topics}

There are so many interesting aspects of $0/1$-polytopes,
and so little time and space. In this section, I am therefore
collecting brief notes about three further topics, together with
pointers to the literature that I'd hope you'll follow.


\subsection{Graphs}


General facts about graphs of polytopes apply in the 
$0/1$-context, but there are new phenomena appearing --- the
most tantalizing perhaps being the Mihail-Vazirani conjecture.
But we start with a basic fact that is true for all (bounded, convex)
polytopes, and hence need not be proved in our more special context.

\begin{Theorem}[Balinski \cite{Balinski}; Holt \& Klee \cite{HoltKlee4}]%
\label{th:B-HK} 
{\rm(1)} The graph of every $d$-dimensional polytope is 
vertex $d$-connected;
that is, there are $d$ vertex-disjoint paths between any pair
of vertices.\\
{\rm(2)} For any generic linear objective function (such that no two 
vertices get the same value), there are $d$ monotone vertex-disjoint paths
from minimum to maximum.
\end{Theorem}

In a setting of general (convex, bounded) polytopes the
first part of this, ``Balinski's Theorem,'' is a classic.
The second part is a rather recent strengthening observed
by Holt \& Klee \cite{HoltKlee4}: it implies the first
part  since for any two distinct vertices of a polytope
we may assume that they are the unique minimal and the 
unique maximal vertex for a generic linear function, 
{\em after a projective transformation\/} \cite[p.~74]{Zpoly}.
One peculiar phenomenon is that this reduction does not
work in a setting of $0/1$-polytopes: 
projective transformations do not preserve $0/1$-polytopes.


The second result for this section is an example of an
important and still unsolved problem from the 
theory of general polytopes (see \cite[Sect.~3.3]{Zpoly})
which becomes quite trivial when specialized to $0/1$-polytopes
--- as was first noticed by Denis Naddef.

\begin{Theorem}[The Hirsch conjecture for $0/1$-polytopes: Naddef~\cite{Na}]%
The diameter of the graph of a $d$-dimensional $0/1$-polytope
$P\sse\R^n$ is at most
\[
\diam(G(P))\ \ \le\ \ d,
\]
with equality if and only if $P$ is (affinely equivalent to) a
$d$-dimensional $0/1$-cube. In particular,
this implies that
\[
\diam(G(P))\ \ \le\ \ n-d,
\]
where $n$ is the number of facets of $P$.
\end{Theorem}

\proof 
We get the first inequality by 
induction on dimension, the case $d=1$ being trivial.
If the two vertices in question lie in a common facet of~$[0,1]^d$,
then we can restrict to the corresponding trivial face of~$P$ of dimension
at most~$d-1$, and we are thus done by induction.
Hence we may assume that $\vv$ and $\uu$ are opposite vertices of~$[0,1]^d$,
and by symmetry only need to consider the case where
$\vv=\zero$ and $\uu=\one$.

But the vertex~$\uu=\one$ is connected to some neighboring vertex $\uu'$,
and this neighbor is contained in some trivial face $P_i^0$,
whose diameter is at most~$d-1$ by induction. Thus
\[
d(\vv,\uu)\ \ \le\ \ d(\vv,\uu')+d(\uu',\uu)\le (d-1)+1=d.
\]
For the second statement, we may assume (using induction on dimension)
that the two vertices in question do not lie on a common facet.
Thus the polytope has at least $n\ge2d$ distinct facets.
%
\proofend

Our third item in this section is a conjecture that's just plain wrong for
general polytopes, but may be true in the $0/1$-setting.

\begin{Conjecture} [Mihail-Vazirani {\cite[Sect.~7]{FederMihail}}]
The graph of every $0/1$-polytope is a good expander. Specifically, for
every partition $V=S\uplus\overline S$ of the vertex set, the polytope
$P(V)$ has at least
\[
E(S,\overline S)\ \ge\ 
\min\{|S|,|\overline S|\}
\]
edges between $S$ and $\overline S$.
\end{Conjecture}

Remark: This may be very false. It does not seem to be trivial.


\subsection{Triangulations}

A very basic question is the following:
{\em How many simplices are needed to triangulate
the $d$-dimensional $0/1$-cube?}
Here the exact answer depends on the exact definitions: for example,
let us assume that we want proper triangulations where all
simplices are required to fit together face-to-face, and not
only subdivisions, or (even worse) coverings.
Let us also assume that we only admit triangulations 
without new vertices. (In general polytopes, new vertices
{\em do\/} help --- see Below et al.\ \cite{BBDLRG}.)
\newcommand{\triang}{\operatorname{triang}}

In this setting, let $\triang(d)$ be the smallest number of
simplices in a triangulation of $[0,1]^d$.
Then we can draw up a little table,
\[
\begin{array}{r|r}
d & \triang(d)\\
\hline
\\[-4mm]
1 &               1 \\
2 &               2 \\
3 &               5 \\
4 &              16 \\
5 &              67 \\
6 &             308 \\
7 &            1493 
\end{array}
\]
combining many earlier results with those of Hughes~\cite{Hugh,Hugh2} 
and 
Hughes \& Anderson \cite{HuAn,HughesAnderson}.
For $d=8$, all we have seem to be the bounds
$5522\le\triang(8)\le11944$.

One of several curious effects in this context is that
{\em not every} $d$-dimensional
$0/1$-polytope can be triangulated into at most
$\triang(d)$ simplices: for example, for the $6$-simensional
half-cube {\tt HC:7-64.poly}, the convex hull of all
$0/1$-vectors of even weight, one knows that the minimal
number of simplices in a triangulation is
$1756 > 1493 = \triang(7)$ (Hughes \& Anderson~\cite{HughesAnderson}).

A lower bound is certainly given by the
maximal volume of a $0/1$-simplex,
\[
\triang(d)\ \ \ge\ \ \frac{d!}{\rho_d},
\]
but this bound is not very good. (For example, for $d=3$ it yields
only $\triang(d)\ge3$). However,
it can be refined by giving greater weight to
the simplices ``near the boundary,'' which have lower volume, but
are needed to fill the $0/1$-cube. A very elegant and powerful
version of such a lower bound was given by Smith~\cite{Smith} 
using hyperbolic geometry.


A good quantity to consider is
\[
\sqrt[{\scriptstyle d}]{\frac{\#\textrm{simplices}}{d!}},
\]
called the {\em efficiency\/} of a triangulation.
This number is at most $1$ for any triangulation that uses
no ``extra vertices.'' Haiman~\cite{Haiman}
showed that the limit
\[
L\ \ :=\ \ \lim_{d\rightarrow\infty}
\sqrt[{\scriptstyle d}]{\frac{\triang(d)}{d!}} 
\]
exists, and that the efficiency of any example can also be
achieved asymptotically, that is,
\[
\sqrt[{\scriptstyle d}]{\frac{\#\textrm{simplices}}{d!}}\ \ \le\ \ L
\]
holds for every triangulation without new vertices.
The best upper bound on~$L$ up to now seems to be the one provided
by Santos \cite{Santos}:
\[
L 
\ \ \le \ \ \sqrt[{\scriptstyle 3}]{\frac7{12}}\ \approx\ 0.836.
\]
One would, however, expect that the limit $L$ is~zero.


\subsection{Chv\'atal-Gomory ranks}

Interesting questions are related to the rounding procedures
of integer programming that try to recover the convex hull
$P_I:=\conv(P\cap\Z^d)$ from an inequality description
of a polytope $P\sse[0,1]^d$.

In particular, Chv\'atal-Gomory rounding steps 
replace $P$ by 
\[
P'\ \ := \ \ \bigcap_{H\supset P} H_I,
\]
where the intersection is taken over all closed halfspaces $H$ that
contain~$P$. The integer closure $H_I$ of a halfspace $H$
is easy to compute:
make the left-hand side of the inequalities integral with
greatest common divisor one, and then round the right-hand side.
It was proved by Chv\'atal that a finite number of such
closure operations lead from a bounded polytope $P$ to its
integer hull --- {\em but how many steps are needed?}
This quantity is known as the {\em Chv\'atal-Gomory rank}
or {\em CG-rank} of the polytope~$P$.
We refer to the thorough treatment by Schrijver~\cite{Schrijver}
for details and references.
\newcommand{\CGr}{\operatorname{CGr}}

Bockmayr, Eisenbrand, Hartmann \& Schulz~\cite{BEHS}
noticed recently that for polytopes in the
$0/1$-cube, $P\sse [0,1]^d$ the Chv\'atal-Gomory rank
is bounded by a polynomial in~$d$.
An improvement of Eisenbrand \& Schulz \cite{EisenbrandSchulz}
establishes that for $P\sse[0,1]^d$ the CG-rank is bounded by
\[
(1+\varepsilon)d\ \ \le\ \ \CGr(d)\ \ \le \ \ 3d^2\log(d)
\]
for some $\varepsilon>0$.

But how about a good lower bound?
Riedel~\cite{Riedel} has implemented a procedure to compute the
CG-rank for polytopes, and he has provided 
explicit, low-dimensional examples $P\sse[0,1]^d$ for which
the CG-rank exceeds the dimension; so we know
\begin{eqnarray*}
\CGr(3)& = & 3\\
\CGr(4)& = & 5\\
\CGr(5)& = & 6 ~~(?)\\
\CGr(6)&\ge& 8\\
\CGr(7)&\ge& 9
\end{eqnarray*}
But can anyone provide
a lower bound that is more than simply linear?


\section{Problems and Exercises}

\begin{enumerate}


\item
Is it true that every simple $0/1$-polytope is a product of simplices?\\
(This question was answered by Kaibel \& Wolff~\cite{KW}.)


\item
Estimate the maximal vertex degree of a $d$-dimensional
$0/1$-polytope.\\ (Hint: {\tt OA:5-18.poly})

\item
Classify the $0/1$-polytopes of diameter $\sqrt{2}$.

\item*%
Bound the maximal number of vertices for a $d$-\allowbreak dimensional
$2$-neighborly $0/1$-\allowbreak poly\-tope.
(Corollary~\ref{cor:neigh} yields an exponential lower bound.)

\item
Show that every $0/1$-polytope without a triangle face
is a $d$-cube.\\
(Volker Kaibel noticed that this follows from a result of
Blind \& Blind \cite{BlindBlind}.)


\item*%
Is it true that a simplicial $0/1$-polytope of dimension $d$
has at most $2^d$ facets?\\
Is it true that every simplicial $0/1$-polytope of dimension $d$
with $2d$ vertices is centrally symmetric and thus is
a cross polytope with exactly $2^d$ vertices?\\
(This is true for $d\le6$, according to Aichholzer's enumerations.)


\item
Estimate the probability that the determinant of a random
$(n\times n)$-matrix with entries in $\Z_2$ vanishes (for large $n$).
Compare your result with that claimed in~\cite{Mu}.

\item
Show that for every fixed $\varepsilon>0$, 
all the trivial faces of a random $0/1$-polytope 
with $(2-\varepsilon)d$ vertices are simplices,
with probability tending to~$1$ for $d\rightarrow\infty$.
\\~~\hfill (Volker Kaibel)

\item\label{Ex:det}
Prove the Szekeres-Tur\'an theorem:
The expected value of the determinant $\det(C)$ of a random
$\pm1$-matrix $C\in\{-1,+1\}^{n\times n}$ is zero, but
the expected value of the squared determinant is exactly~$n!$:
\[
E(\det(C)^2)\ \ = \ \ n!.
\]
Hint, by Bernd G\"artner: Use 
$\det (C) = \sum_{i=1}^{n} (-1)^{i-1}c_{1i}\det(C_{1i})$, 
and analyze the expected values of the summands in 
\[
\det (C)^2\ \ =\ \ \sum_{i=1}^{n} (\det(C_{1i}))^2
\ + \  \sum_{i\neq j} (-1)^{i+j}c_{1i}c_{1j}\det(C_{1i})\det(C_{1j}).
\]

\item
What is the largest absolute value of the determinant of an
$n\times n$ matrix with coefficients in $\{-1,0,1\}$?
With coefficients in the interval $[0,1]$? With coefficients in~$[-1,1]$?
\\
(It is reported that this is a question that was asked by L. Collatz
at an international conference in 1961, and answered a year later by
Ehlich \& Zeller \cite{EZ}. Your answer should be in terms of $\rho_n$
resp.\ $\rho_{n-1}$.) 

\item
Show that $\CUT(k)$ is ($0/1$-isomorphic to) a face
of $\CUT(n)$, for $k\le n$.

\item
Prove that $[\zero,\one]$ is an edge of the correlation polytope $\COR(n)$.

\item \label{Ex:6-simplicial}
Show that $\CUT(n)$ is \emph{$6$-simplicial}: every
$5$-dimensional face is a simplex.

\item
Show that the \emph{metric polytope} 
\[
\MET(n):=
\left\{
\begin{array}{ll}
X\in[0,1]^{d}: & x_{ij}-x_{ik}-x_{jk}\le 0\textrm{~~ and}\\
               & x_{ij}+x_{ik}+x_{jk}\le 2\textrm{~~ for distinct
                         } i,j,k\in[n]
\end{array}
\right\}
\]
is an $LP$-relaxation of $\CUT(n)$: it satisfies
$\CUT(n)=\conv(\MET(n)\cap\Z^{d})$,
where $d=\binom{n}{2}$ is the dimension of 
$\CUT(n)\subseteq\MET(n)\sse\R^d$.
\\
*Estimate the CG-rank of $\MET(n)$.

\item
How do the inequalities for a $0/1$-polytope transform into the
inequalities for the corresponding $({+}1/{-}1)$-polytope?

\item
Give more and better examples of ``large'' coefficients appearing in the
facet-defining inequalities of $0/1$-polytopes.



\item
Show that every triangulation of $\Delta_k\times\Delta_{\ell}$ 
without new vertices has exactly $\binom{k+\ell}{k}$ facets.

\item*\label{prob:reg_sx}%
For which dimensions $d>1$ and integers~$k$ ($1\le k\le d$)
does there exist a regular $d$-dimensional $0/1$-simplex
of edge length~$\sqrt{k}$?\\
(Show that this is equivalent to the existence of 
a matrix $M\in\{0,1\}^{d\times d}$ with 
$M^tM=\frac{k}{2}(I_d+\one^t\one)$, so in particular
$k$ must be even. Show that the case $k=d$ is equivalent
to the famous Hadamard determinant problem.)

\item*%
For which $d$ is there a regular $d$-dimensional $0/1$-cross polytope?

\item
For which $E\subseteq\binom{[n]}{2}$ is
$P(E)=\conv\{\ee_i+\ee_j:\{i,j\}\in E\}$ a simplex?\\
Show that every such simplex of dimension $d=\binom{n}2$
has normalized volume $\frac{2^k}{d!}$ for
some $k\ge 0$. \cite{DLST}

\item
Estimate the volumes of the Birkhoff polytopes
\[
B_{n+1}:=\{X\in[0,1]^{n\times n}:
\one^t X\le\one^t,\ 
X \one\le\one,\ 
\one^t X\one\ge n-1\}.
\] 
(see {\tt BIR3:4-6.poly}, {\tt BIR4:9-24.poly}, \ldots).
The exact value of the volume of~$B_{n+1}$, which is some integer
divided by $n^2!$, is known for~$n\le7$, due to
Chan \& Robbins~\cite{ChanRobbins}.


\end{enumerate}

\section*{Acknowledgements}

Thanks to Bernd G\"artner, Carsten Jackisch, Fritz Eisenbrand, 
Francisco Santos, Gerald Stein, Imre B\'ar\'any, 
Ji{\v r}{\'\i} Matou\v{s}ek, Marc Pfetsch, Michael Joswig, Noga Alon, 
Oswin Aichholzer, Thomas Voigt and Volker Kaibel
for so many helpful discussions and useful comments.

\newcommand{\reference}[4]{\bibitem{#1}{\sc #2 \it #3 \rm #4}}
\newcommand{\inproc}[1]{{\rm in:} ``{\rm #1}''}

\begin{small}

\end{small}
\end{document}